\documentclass{article}
\usepackage{mathptmx}
\usepackage{graphicx}
\usepackage{float}
\usepackage[tight]{subfigure}
\usepackage{amsmath, amssymb}

\usepackage{parskip}

\def\bw{\mbox{\boldmath$\omega$}}
\def\bx{\mathbf{x}}

\def\bD{\mathbf{D}}
\def\bn{\mathbf{n}}
\def\bg{\mathbf{g}}

\def\bpsi{\mbox{\boldmath$\psi$}}

\def\bu{\mathbf{u}}
\def\bU{\mathbf{U}}
\def\bnabla{\mbox{\boldmath$\nabla$}}

\title{A Vortex Method for Bi-phasic Fluids\\ Interacting with Rigid Bodies}

\author{Mathieu Coquerelle\thanks{e-mail: Mathieu.Coquerelle@imag.fr }\\ GRAVIR/IMAG and LMC/IMAG 
\and J\'er\'emie Allard\thanks{e-mail: Jeremie.Allard@imag.fr }\\ GRAVIR/IMAG %
\and Georges-Henri Cottet\thanks{e-mail: Georges-Henri.Cottet@imag.fr }\\ LMC/IMAG %
\and Marie-Paule Cani\thanks{e-mail: Marie-Paule.Cani@imag.fr }\\ GRAVIR/IMAG
}
\date{}

\begin{document}

\maketitle


\begin{abstract}
We present an accurate Lagrangian method based on vortex particles,
level-sets, and immersed boundary methods, for animating the interplay
between two fluids and rigid solids.  We show that a vortex method is
a good choice for simulating bi-phase flow, such as liquid and gas,
with a good level of realism.  Vortex particles are localized at the
interfaces between the two fluids and within the regions of high
turbulence.  We gain local precision and efficiency from the stable
advection permitted by the vorticity formulation. Moreover, our
numerical method straightforwardly solves the two-way coupling problem
between the fluids and animated rigid solids.  This new approach is
validated through numerical comparisons with reference experiments
from the computational fluid community. We also show that the visually
appealing results obtained in the CG community can be reproduced with
increased efficiency and an easier implementation.

\end{abstract}

\section{Introduction}

Fluids interacting with solid objects are a common, yet fascinating every-day life experience. Our tendency to stare at turbulent liquids and at smoke dynamic behavior, or to observe different objects splashing into water makes us very critical when we see such phenomena reproduced with a computer. Yet, the demand for plausible fluids in computer-generated movies and games is high.
This has led many computer graphics researchers to tackle the challenge, although still considered as one of the most difficult problem in the computational fluids community. 

The interactions we generally observe in real life involve several components: typically, water with a free surface in contact with the air, which interacts with both still and moving rigid bodies.
These interactions result in really complex phenomena (turbulences, splashes, bubbles, interesting subsequent motion of the solids), due to the interplay between the two fluids and the solids.
This paper presents a Lagrangian approach, based on vortex particles, for accurately, yet efficiently simulating this interplay. We show that vortex particles
provide a good solution to the modelling of bi-phase flow with a liquid component.
Since they concentrate the computational power in discontinuous and turbulent regions, the complex phenomena occurring at the interface between the two fluids are represented with a good level of precision.
Meanwhile, the numerical method we present straightforwardly solves the two-way coupling problem between the fluids and animated rigid solids, without the need for complex boundary conditions.

\subsection*{Previous Work}

The simulation of fluids and of their interplay with solid objects has attracted a lot of attention in the CG community within the past few years. Great advances were recently made towards this goal.

\cite{foster96realistic} presented the first 3D simulation of liquids.
They relied on an Eulerian formulation to solve the Navier-Stokes equations for incompressible fluids.
A semi-Lagragian method for advecting the fluid which guaranties the stability of the simulation was then proposed by \cite{stam99stable}.
Because of high numerical dissipation, these kind of simulations loose part of their vorticity over time, a problem which was tackled in \cite{fedkiw01visual}.
\\
Lagrangian methods based on particles provide a good alternative to Eulerian simulations, since they enable to dynamically follow the fluid. The \textit{Smoothed Particles Hydrodynamics} (SPH) formulation was adapted by \cite{desbrun98active} in order to transport an implicit boundary with surface tension.
More recently, \cite{muller05particle} used this formulation to simulate fluid-fluid interactions such as boiling water and managed to obtain interesting phenomena such as air bubbles in water.
Closer to the Navier-Stokes equations for incompressible fluids, \cite{premoze03particlebased} proposed a particle based solution called \textit{Moving Particle Semi-Implicit} (MPS).
While straightforwardly able to track multiple fluids, particles-only methods mostly suffer from the difficulty to conserve the dynamic properties of the fluids, especially when dealing with boundaries.
\\
Among these approaches, vortex methods~\cite{gamito95two,Park_2005_5216,AN05} were used in for animating gaseous phenomena. They are particularly interesting since they focus the resolution in the regions of high turbulence. They rely on the vorticity formulation of the Navier-Stokes equations, well known in applied mathematics for being an accurate alternative to Eulerian methods (see \cite{C-CotKou00}).
Vortex particles were also introduced to add subgrid turbulences on top of an Eulerian simulation of liquids or gazes ~\cite{selle05vortex}, impressively counteracting the numerical dissipation due to the underlying Eulerian solver. Up to now, no technique was developed in graphics to simulate liquids or multi-phases fluids from the vorticity formulation of Navier-Stokes. This may be due to the inherent cost of particle methods (from $O(N^2)$ to $O(N log N)$ using \cite{lindsay}), more expensive than Eulerian methods since a large number of particles ($N\gg1000$) are necessary to reach a good level of precision. One of the contributions of this paper is to show that when adequately implemented, the vortex formulation achieves accurate yet efficient simulations. 
\\

Animation and visualization of water in contact with the air require the precise transport and rendering of the interface between them.
Because of the discontinuity of the fluids' physical properties at the interface, the resolution of the Navier-Stokes equations was most of the time restricted to the liquid component and to the interface region~\cite{enright02animation,enright05fast}. However, this prevents simulating some interesting phenomena such as air bubbles inside the liquid.
A multi-phase method was recently presented to take care of both fluids, achieving impressive bubbles movements \cite{hong05discontinuous}.
As in this last work, we simulate bi-phase fluids. We present a different solution, based on the vorticity formulation, which solves the problem in an intuitive and efficient way.
\\

Fluid-structure interaction is considered as a really tough problem by both the CG and applied mathematics communities.
The most difficult issue is to compute the fluid's velocity at the objects boundaries while ensuring that no fluid penetrates any obstacle.
\cite{genevaux03simulating} presented one of the first methods in CG for generating interaction forces between fluids and rigid solids.
Recently, \cite{carlson04rigid} provided a solution that treats the rigid body as a fluid and extracts its associated movement from the fluid's velocity.
\cite{guendelman05coupling} presented the first technique able to handle the interactions with thin deformable and rigid objects represented by meshes.
This method was further improved by \cite{losasso05melting} to allow the melting and burning of these solids.
Closer to our vortex particles approach, \cite{Park_2005_5216} used the panel method to impose no-slip and no-through constraints on the fluid by emitting new particles at the solid's boundaries.
This one-way interaction method, while well adapted for gases, is based on an explicit representation of the object which, if the density of particles was too low, does not prevent the flow from penetrating.
Although achieving impressive results, these methods are sometimes hard to implement or require complex treatments where the object touches the fluid. Furthermore, they do not provide any validation of the correctness of the fluid's behaviors at the boundaries.
Although in the spirit of the 'rigid fluid' method \cite{carlson04rigid}, the solution we present avoids the explicit tracking of the fluid/body interface.
Instead, we use a level set formulation on top of a vorticity creation algorithm to apply forces and account for the continuity of velocity at the fluid/body interface.

\subsection*{Overview}

Our motivation for using vortex particles lies in the following features: they allow to concentrate numerical efforts on regions of interest and remain stable for large time-step values.
Their main limitation, the associated computational cost, can be alleviated by using an auxiliary 3D grid. This is done while keeping the vorticity formulation of the equations and remaining physically accurate.

Our first contribution is to use vortex particles to simulate bi-phase fluids such as liquid and gas. As we show, the vortex formulation is very well adapted to tackle this problem, since the vortex particles are created near the interface between the two fluids.
The second contribution is a novel algorithm to simulate the two-way interactions of the bi-phase fluid with animated rigid bodies.
Our approach provides a physically sound and clear-cut fluid-solid model thanks to an algorithm that is both robust and easy to implement.

Section \ref{sect:vortex particles} reviews vortex methods.
Section \ref{sect:bi-phasic} explains how we use them to simulate bi-phase flow, such as water in contact with air.
Section \ref{sect:interactions} focuses on the animation of solids interacting with the fluids.
Section \ref{sect:validation} validates the method through numerical comparisons with reference experiments from the computational fluid community.
Lastly, Section \ref{sect:results} shows that the visually appealing results obtained in the CG community can be reproduced with an increased efficiency.

\section{Vortex Particle Methods}
\label{sect:vortex particles}

In most fluid simulations, only a part of the flow has an interesting behavior.
The vorticity formulation derived from the Navier-Stokes equations allows to 
focus the computational cost on the region of interest using vortex particles.
\\
This section gives a quick overview of the vortex particles method, we refer 
the reader to \cite{C-CotKou00} and the references therein for detailed 
numerical analysis and discussions of this method.

\subsection*{Definition}

We start from the incompressible viscous 3D Navier-Stokes equations:
\begin{eqnarray}
\label{eq:ns3D}
\bu_t + (\bu\cdot \bnabla)\bu -\nu\Delta \bu + \nabla p = 0
\\
\label{eq:ns3Dincomp}
\nabla \cdot \bu = 0
\end{eqnarray}
where $\bu_t$ is the velocity's time derivative of the fluid's velocity field $\bu$, $\nabla\cdot\bu$ 
is the divergent of $\bu$, $p$ is the pressure and $\nu$ the kinematic 
viscosity of the fluid. 
The vorticity is defined as the curl of the velocity:$$\bw=\nabla\times\bu$$ 
(note that in the particular case of a 2D fluid, the vorticity is a scalar).
Taking the curl of equations (\ref{eq:ns3D}) and (\ref{eq:ns3Dincomp}) (after 
some linear algebra) leads to the vorticity formulation of the Navier-Stokes 
equations:
\begin{equation}
\label{ns}
 \bw_t + (\bu\cdot \bnabla)\bw = (\bw\cdot \bnabla)\bu +\nu\Delta \bw.
\end{equation}
Solving this single equation is equivalent to solving equations (\ref{eq:ns3D}) 
and (\ref{eq:ns3Dincomp}).
The term $(\bu\cdot \bnabla)\bw$ represents the transport of the vorticity by 
the fluid's velocity, the left-most term of the right-hand side of equation 
(\ref{ns}) represents the stretching (change of orientation) of the vorticity 
vector while the latter term represents the vorticity diffusion due to 
viscosity.
\\
While it would be possible to solve this equation on a traditional 3D grid, 
the use of particles to transport vorticity permits to gain in precision and 
efficiency.
Contrary to velocity which is mostly non zero in the whole domain, vorticity 
is localized in region of turbulences even if there can be high velocities 
everywhere.
In consequence, by the use of particles, computational precision is focused 
where vorticity exists.
Another advantage is that the advection of vorticity is not subjected to the 
time-step constraint inherent to Eulerian methods, which stands that $\delta t 
< |\bu|_{max}$.
\\
Vortex particle methods consist in representing the vector field $\bw$ by a 
set of particles:
$$\bw(\bx)=\sum_p v_p \bw_p \zeta(\bx-\bx_p)$$
where $\bx_p$, $\bw_p$ and $v_p$ are respectively the location, strength and 
volume of particle $p$ and $\zeta$ is a smooth distribution function, 
typically a Gaussian.
Due to the incompressibility constraint, the volumes $v_p$ remain constant.
Rewriting equation (\ref{ns}) in a Lagrangian formulation, the particles' location and 
strength are integrated using:
\begin{eqnarray}
\label{eq:partvelox}
\frac{\bD{\bx}_p}{\bD t}=\bu(\bx_p,t),
\\
\label{eq:partstren}
\frac{\bD{\bw}_p}{\bD t} = (\bw\cdot \bnabla)\bu +\nu\Delta \bw,
\end{eqnarray}
where $\bD q/\bD t=\partial q / \partial t + (\bu\cdot\nabla)q$ denotes the 
rate of change of a quantity $q$ in the Lagrangian frame of a fluid element 
advected by the fluid.
Particles' velocity and vorticity derivatives have to be determined in a self-
consistent way from the vorticity field.

\subsection*{Coupling Vortex Particles with a Grid}
As we plan to simulate fluids which may become turbulent and thus require many 
particles, we use both vortex particles and a 3D grid: firstly, 
spatial differentiations are cheaper on a grid; secondly it guaranties an 
approximately constant cost when solving the equations; lastly, it solves with 
no extra cost the problem of redistributing the particles over time.
This last point is important because particles advected by the fluid will 
naturally tend to cluster in some area or to move away from each other thus 
leaving unresolved spaces between them.
\\
A fast and accurate solution is to superpose a uniform grid on the particle 
distribution: at each time step particles are remeshed on that grid (thus, the 
particles' volume $v_p=h^3$ where $h$ is the grid's cells spacing).
Remeshing at every time-step can be interpreted as a class of high order finite-difference schemes, as long as it is done with high order redistribution schemes (see \cite{chaniotis} for example).
\\
Given the vorticity field, velocity is computed on the grid by means of a fast 
grid solver (in our experiments we used the FishPack library \cite{fishpack} in 
order to solve the following Poisson equation:
\begin{equation}
\label{eq:streamfunction}
\Delta \bpsi = - \bw,
\end{equation}
to obtain the so-called stream function $\bpsi$ which is differentiated to 
obtain the velocities :
\begin{equation}
\label{eq:stream2velox}
\bu = \nabla \times \bpsi.
\end{equation}
This solution to compute the velocity field from the vorticity is 
faster than the particle based solution used in \cite{Park_2005_5216} when the number of particles increases.\\
In practice particles are advanced with a second or fourth order Runge-Kutta 
method.
During each substep velocity and vorticity's time derivative (the right hand 
side of the vorticity equation (\ref{ns})) are interpolated from the grid onto 
the particle locations.
The accuracy of the overall algorithm heavily relies on the quality of the 
interpolation formulas used to remesh particles on the grid and to interpolate 
back fields at particle locations.
It is common practice in CFD to use smooth interpolation formulas which 
preserve moments of order up to 2.
The stencil on which particles are remeshed extends to the four nearest grid 
points in each direction.
The cost of the algorithm thus scales linearly with the number of particles. A number of numerical validations on challenging test cases in CFD has allowed to attest the accuracy of remeshed particle methods. These studies are backed by  a recent theoretical analysis which demonstrate that remeshed particle methods are equivalent to a class of high order finite-difference schemes not subject to CFL conditions (\cite{CotWey06}).  

\section{Bi-phase fluids}
\label{sect:bi-phasic}

We present here an intuitive method for simulating bi-phase fluids which 
benefits from the vortex particles formulation.
Before explaining the particularity of bi-phase fluids, we first consider the 
vorticity problem in the presence of a fluid of variable density.

For variable density and viscosity flows, the vorticity equation (\ref{ns}) 
must be replaced by:
\begin{align}
\label{eq:nsvd}
\bw_t + (\bu\cdot \bnabla)\bw = (\bw\cdot \bnabla)\bu & + \nabla\cdot
(\nu\bnabla \bw) + {(\nabla p \times \nabla \rho) \rho^{-2}}
\notag	
\\
&+ \bnabla\times (\rho \bg) + \nabla\cdot[(\bnabla\times\nu)\bnabla\bu]
\end{align}
 where $\bg$ is the gravity field, $\rho$ is the density and $\nu=\nu(\rho)$ 
is the variable kinematic viscosity.
This equation has to be supplemented by a transport equation for $\rho$ or for 
$\nabla \rho$.
\\
This system is often simplified by assuming small density variations: 
in the so-called Boussinesq approximation, pressure and velocity terms in (\ref{eq:nsvd}) disappear and we are left with 
\begin{equation}
\label{boussinesq}
\bw_t + (\bu\cdot \bnabla)\bw = (\bw\cdot \bnabla)\bu + \nabla\cdot(\nu\bnabla 
\bw) + \bnabla\times (\rho \bg).
\end{equation}
It is important to note that, in the case of two fluids with constant 
viscosity and variable density, this equation shows that vorticity is produced 
where density gradients are located, that is at the interface between the 
fluids.
Thus computation is localized in this narrow band which clearly reduces the 
theoretical cost compared to traditional methods.
Accordingly to (\ref{boussinesq}), the particle's vorticity change in equation (\ref{eq:partstren}) becomes:
\begin{equation}
\label{eq:partstrengrav}
\frac{\bD{\bw}_p}{\bD t} = (\bw\cdot \bnabla)\bu + \nabla\cdot(\nu\bnabla 
\bw) + \bnabla\times (\rho \bg).
\end{equation}
We now consider the simpler problem of a bi-phase fluid in a domain $\Omega$ 
where the density (resp. viscosity) takes only two different values: $\rho_1$ 
(resp. $\nu_1$) in a domain $\Omega_1$ and $\rho_2$ (resp. $\nu_2$) in a 
domain $\Omega_2$ ($\Omega=\Omega_1 \cup \Omega_2$).
We use a level set (noted $\phi$) to capture the interface between those two 
domains:
$$
\phi(\bx) \left\{
\begin{array}{l}
< 0 \mbox{ for $\bx \in \Omega_1$ },
\\
= 0 \mbox{ for $\bx \in \Omega_1 \cap \Omega_2$ },
\\

> 0 \mbox{ for $\bx \in \Omega_2$ }.

\end{array}
\right.
$$
Typically, we initialize $\phi$ as a signed distance.
This level set function satisfies a transport equation which can be solved 
either in its primitive form:
  \begin{equation}
\label{ lsf1}
\phi_t+ (\bu\cdot\bnabla) \phi =0
\end{equation}
or in its gradient (vector) form:
\begin{equation}
\label{lsf2 }
(\bnabla \phi)_t+ (\bu\cdot\bnabla) \bnabla\phi = -(\bnabla\phi\cdot\bnabla)\bu
\end{equation}
The latter equation is very similar to the vorticity equation and thus gradient of $\phi$ and $\bw$ are both located near the interface.

In this case particles carry strengths of vorticity and of $\nabla\phi$, from which we can deduce $\nabla\times\rho$ needed in equation (\ref{eq:partstrengrav}) (see below).

But contrary to \cite{B-Cot03}, we advect the level set $\phi$ on the grid with a semi-Lagrangian method.
The surface tension is defined in terms of the curvature of $\phi$ (see~\cite{hou}):
$$
\tau \kappa \zeta(\phi) \nabla \phi.
$$
where $\tau$ is the surface tension coefficient.
This term is added to the vorticity equation (\ref{eq:partstrengrav}).
The level set is interpolated onto the particles during the remeshing step.
\\
In order to be able to solve equation (\ref{eq:partstrengrav}), we need to know the density and viscosity fields.
They are defined in the whole domain thanks to a smoothed version $H$ of the Heaviside function (with values $0$ and $1$ respectively 
for positive and negative arugments, and where $\varepsilon$ is a small 
parameter of the order of the grid-size) based on the level set:
$$
\rho = \rho_1 H(\phi/\varepsilon) + \rho_2 (1-H(\phi/\varepsilon)).
$$
The same convention is used for viscosity.
From this formulation, it is trivial to deduce the density gradient $\nabla\rho$: $\nabla\rho = (\rho_1-\rho_2)\nabla\phi\zeta(\phi/\varepsilon)/\varepsilon$.
Obtaining $\nabla\times\rho$ is straightforward.
\\
As the computation of \ref{eq:nsvd} is much more involved and expensive then equation \ref{boussinesq}, we have implemented the latter.
Despite the fact that the Boussinesq approximation applies on small density 
variations, we have found that simulating a bi-phase fluid with a high 
difference of density (e.g. $\rho_{air}\approx1$ and $\rho_{water}\approx1000
$) using this technique provides visually realistic results (for both fluid's 
dynamic and interface localization).

\subsection*{Algorithm}
\begin{figure}[t]
\centering
\subfigure[Steps 1 and 2: from vorticity to velocity (on the grid, vorticity 
is denoted by red circles). The magenta curve represents the interface between 
the two fluids.]{
\label{fig:algo_a}
\includegraphics[width=.4\textwidth]{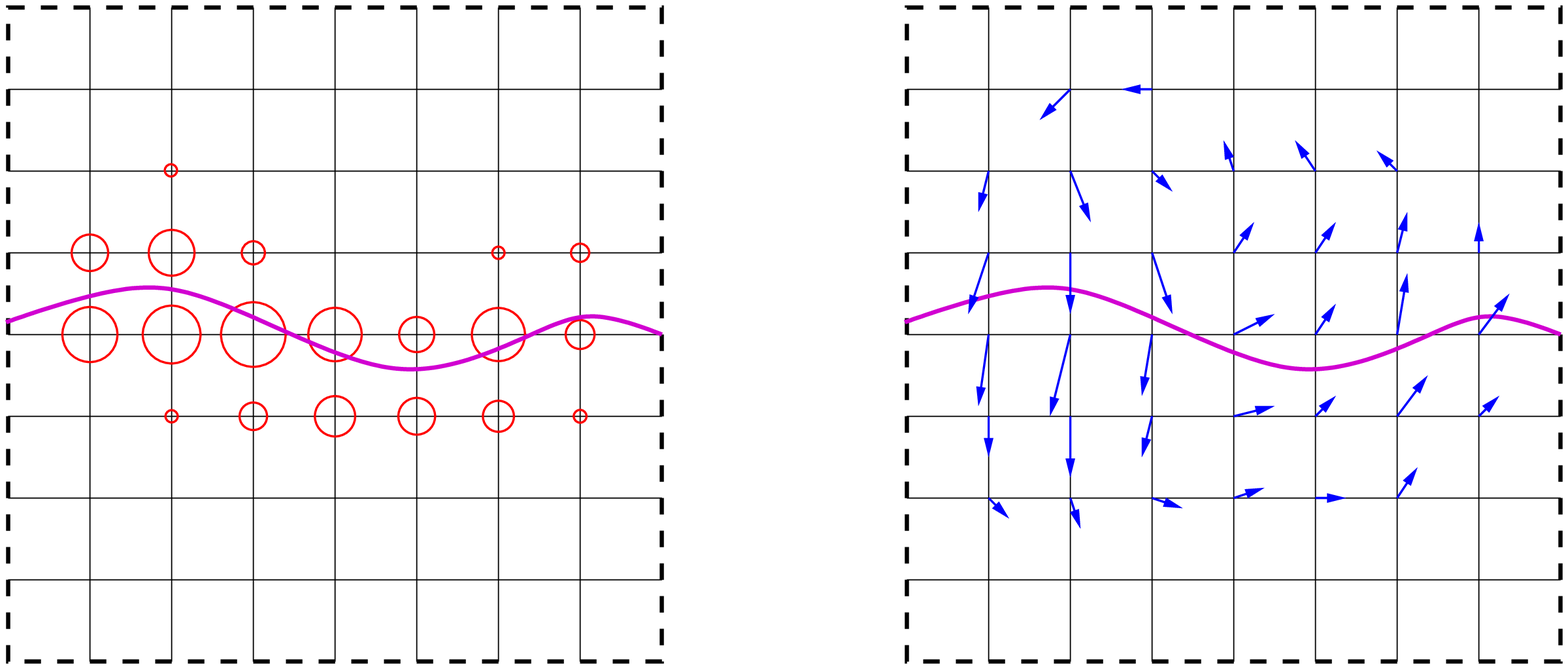}
}\hspace{0.2cm}
\subfigure[Steps 4 and 5: interpolation and advection of particles (particles' 
vorticity is denoted by green discs). The thick dashed (resp. thin dotted) 
magenta represents the advected level set (resp. level set at previous time).]{
\label{fig:algo_b}
\includegraphics[width=.4\textwidth]{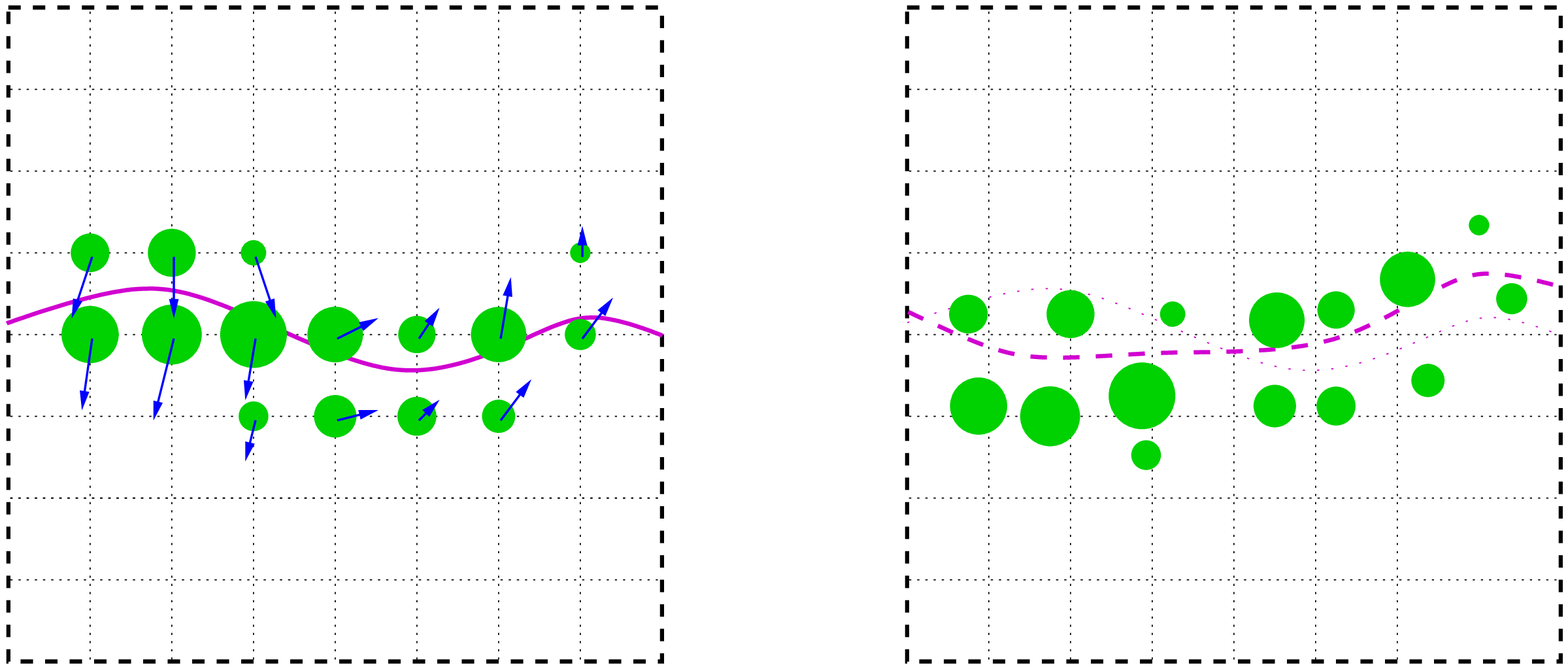}
}
\subfigure[Step 6 (left figure): particles (in transparent green on 
background) are distributed onto the grid. Step 8: remeshing step: create 
fresh particles where vorticity is strong enough (green crosses denote the fact 
that particles were not created because vorticity was too small).]{
\label{fig:algo_c}
\includegraphics[width=.4\textwidth]{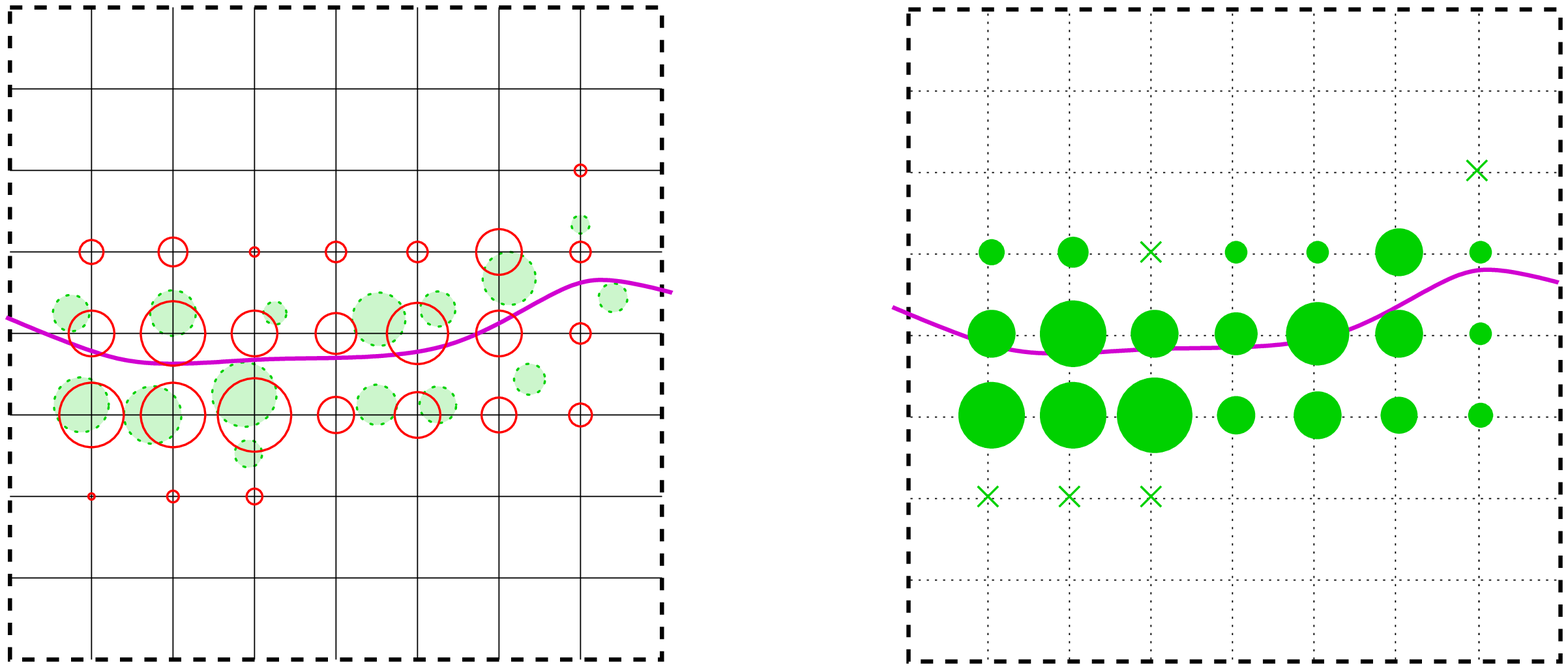}
}
\caption{Algorithm summary}
\label{fig:algo}
\end{figure}

In summary, the following process is applied to advance from time $t_n=n\Delta 
t$ to $t_{n+1}$. Only steps 5 and 8 differ from the standard implementation we presented in Section \ref{sect:vortex particles}:
\begin{enumerate}
\item Find the stream function by solving equation (\ref{eq:streamfunction}) 
on the grid,
\vspace{-0.2cm}
\item Compute the velocity $\bu^{n+1}$ from equation (\ref{eq:stream2velox}),
\vspace{-0.2cm}
\item Compute the stretching term $(\bw \cdot \nabla) \bu$ on the grid,

\vspace{-0.2cm}
\item Interpolate velocity and the stretching term on the particles,
\vspace{-0.2cm}
\item Advance the particles: advect them with velocity (equation (\ref{eq:partvelox})) and vorticity change (equation (\ref{eq:partstrengrav})); Advect the fluid's level set with the velocity on the grid,
\vspace{-0.2cm}
\item Distribute the particles onto the grid to get the advected vorticity,
\vspace{-0.2cm}
\item Compute and integrate the viscosity term $\nabla\cdot( \nu\Delta\bw )$ 
on the grid,
\vspace{-0.2cm}
\item Create fresh vortex particles, carrying $\bw$ and $\nabla\phi$, from the grid if vorticity is greater than 
a threshold.
\end{enumerate}
Note that we used a \textit{viscous splitting algorithm}~\cite{chorin78} to handle separatly the advection and the diffusion of the vorticity.
Please refer to figure \ref{fig:algo} for a visual explanation of the 
algorithm.

\section{Coupling Fluids with Animated Solids}
\label{sect:interactions}

We propose a new approach which, with the adjoint power of the vortex 
particles and level set methods, can be used to compute fluid-solid 
interactions in an intuitive manner and for a low additional cost.
Our approach is to consider the fluid-body system, during the interaction 
processing step, as two fluids of different densities and different 
constitutive laws.

Unlike the method in \cite{carlson04rigid}, our method computing the forces applied by 
the solid on the fluid (and vice-versa) does not require the explicit tracking 
of the interface.
We instead use a level set method combined with a vorticity creation 
penalization term that enforces velocity continuity at the fluid-solid 
interface.
\\

In order to simplify the equations we assume here a single solid, the 
extension to multiple solids being straightforward.

At time zero the body is represented by the zero surface of a level set 
function noted here $\varphi$ in order to not confuse the reader with the 
fluid's level set $\phi$.
Typically, we initialize the level set function as a signed distance to the 
body boundary, negative inside the solid.
This is the only costly part of the method but it does not influence its 
efficiency since it is done as a precomputation.
\\
We now denote by $\widetilde{\bu}^{n+1}$ the velocities found in step 2 of the 
previous algorithm.
After this step, we project the velocities in the grid onto rigid body 
velocity $\bU^{n+1}$ inside the solid using the following formula:
\begin{equation}
\label{rb1} \bU^{n+1}= (\int \widetilde{\bu}^{n+1} \,\, H(\varphi/\varepsilon) 
\,d\bx) / (\int  H(\varphi/\varepsilon) \, d\bx).
\end{equation}
$\bU^{n+1}$ stands for the mean velocity inside the solid.
A similar equation defines the mean vorticity. Integration is performed on the 
whole domain covered by the solid).
At the same time, we enforce continuity of velocities at the fluid-solid 
interface with a penalization technique that we will describe below.
\\

At this stage, velocity and vorticity in the whole domain are obtained by the 
formula :
\begin{equation}
\label{final_velocity}
\bu^{n+1}=\bU^{n+1}  H(\varphi/\varepsilon) + \widetilde{\bu}^{n+1}\,(1- H
(\varphi/\varepsilon)).
\end{equation}
A similar formula is used to obtain the vorticity field $\bw^{n+1}$.
Step 3 and 4 of the algorithm now use these fields for computing the 
stretching term and for the interpolations on the particles.
\\
Step 5 is modified in two ways.
First, we take care of the solid's density exactly in the same way as we have 
done for bi-phase fluids.
Thus, particles now carry the gradient of a fluid with three different 
densities.
Secondly, we need to advect the solid's level set.
This is done just after step 5 of the algorithm of Section \ref{sect:bi-phasic}.
As the solid's velocity gives a translation and a rotation terms, the latter 
problem is easily solved by simply applying the rigid transformation between 
the initial level set position at time $t_0$ and its current position.
In consequence, this scheme does not suffer from temporal diffusion and only 
implies a low diffusion depending on the order of the spatial interpolation.
In our experiments, we use a simple first order interpolant.
\\

It remains to explain the penalization method to enforce velocity continuity.
We use the penalization method of \cite{bruneau}:
assume we want to solve Navier-Stokes equations in a fluid domain outside a 
solid domain $S$, with velocity $\bar \bu$ at the interface.
The penalization model reads
\begin{equation}
\label{penalization_u}
\bu_t + (\bu\cdot \bnabla)\bu -\nu\Delta \bu + \nabla p = \lambda\chi_S
(\bar\bu-\bu),
\end{equation}
where $\chi_S$ denotes the characteristic function of the solid domain $S$ ($1
$ inside and $0$ outside) and $\lambda\gg 1$ is a penalization parameter.
Adapting this method to the vorticity formulation for our problem leads to the 
substitution of the initial vorticity equation (\ref{ns}) by:
\begin{align}
\label{penalization_omega}
\bw_t+ (\bu\cdot \bnabla)\bw = (\bw\cdot\bnabla)\bu + \nu\Delta \bw &+ 
\lambda\chi_S(\bar{\bw}-\bw)
\notag
\\
 &+ \lambda\delta_{\partial S}\,\bn\times(\bar{\bu}-\bu),
\end{align}
where $\delta_{\partial S}$ is the 1D Dirac mass supported by the solid 
boundary and $\bn$ is the unit normal, pointed inward.
Thus, to be consistent with the level set representation for solids, the final 
equation is:
\begin{align}
\label{penalization_ls}
\bw_t+ (\bu\cdot \bnabla)\bw = (\bw\cdot\bnabla)\bu + \nu\Delta \bw  &+ \lambda 
H(\varphi/\varepsilon) (\bar{\bw}-\bw)
\notag
\\
& +\lambda \frac{\zeta({\varphi/\varepsilon})}{\epsilon} \,\bnabla \varphi \times (\bar{\bu}-\bu).
\end{align}

The penalization term tends to cancel the vorticity difference inside the body 
and adds vorticity at the fluid-solid interface.
This term is computed on the grid after the computation of the fluid's velocity at step 2 of the algorithm.
It is then interpolated on the particles in step 4.
We use an explicit scheme to discretize it and choose $\lambda=1/\delta t$ to enforce the stability of this scheme.

\section{Validation}
\label{sect:validation}

In this section, we validate the robustness and the physical accuracy of our method for computing the interactions between fluids and solids.
We compare our results to the reference work from \cite{kern05} obtain with the STAR-CD\footnote{Get more informations on STAR-CD at http://www.cd-adapco.com/sitemap.html} software, a well-known software in the Computational Fluid Dynamic community, implementing fluid-body interactions from \cite{farhat00}.
\\
For this purpose, we first study the 2D case of a falling cylinder (which may be seen as a falling disk) in a fluid of constant density, submited to gravitational forces.
We use the following parameters: boundary conditions are periodic in a square domain of size $1.0\times1.0$, cylinder's radius is $0.1$, fluid's and cylinder's density are respectively $\rho_{fluid}=1$ and $\rho_{cylinder}=2$, kinematic viscosity $\nu$ is $0.001$, gravity is $-1 \mathbf{e}_y$ and $\varepsilon=2\delta x=2\delta y=h$ ($\varepsilon$ is used for computing the smoothed Heaviside function and Dirac distribution values).
Time steps (resp. cells size) are ${\delta t}_{128}=0.01$ (resp. $h_{128}=7.8125\cdot10^{-3}$), ${\delta t}_{256}=0.0038$ (resp. $h_{256}=3.90625\cdot 10^{-3}$) and ${\delta t}_{300}=0.0027$ (resp. $h_{300}=3.33\cdot 10^{-3}$) for grids of dimension $128\times128$, $256\times256$ and $300\times300$ respectively.
The two last time steps are constrained by the diffusion stability condition for an exptheslicit scheme.
\\
Their simulation uses a time step of ${\delta t}_{Kern}=0.0005$ and a 2D radial grid following the cylinder, 15800 cells, $80\times200$ nodes (radial$\times$tangential) and first cell at $\delta r = 0.65\cdot10^{-3}$.
This type of grid is time consuming, but well adapted for a sharp resolution of this particular problem as it permits to focus the computations near the cylinder where precision is the most important.
Both ours and Kern's scheme are $2^{nd}$ order in space and Kern's method is $1^{st}$ order in time while we are using a Runge-Kutta 2 method for advecting the particles.
\\

\begin{figure}[t]
\centering
\includegraphics[width=.65\textwidth]{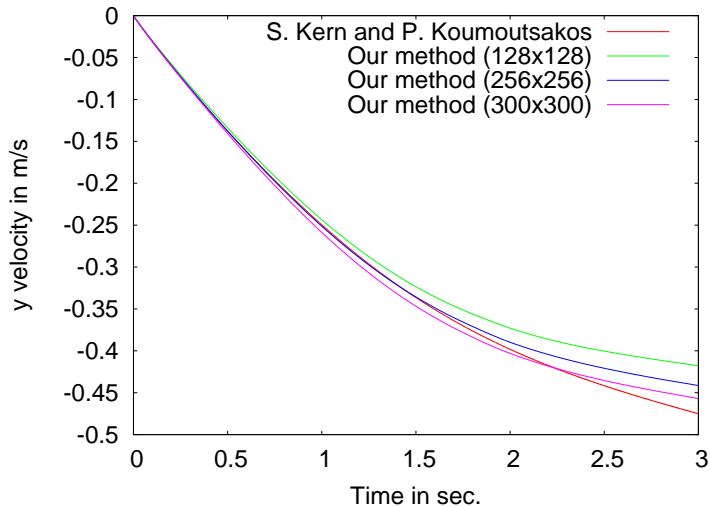}
\caption{$y$ velocities of the cylinder under gravity influence with a RK2 advection for our method, compared to the reference curve.}
\label{fig:velox}
\end{figure}

The $y$ component of the velocity obtained with our method (see figure \ref{fig:velox}) is quite similar to the one obtained with \cite{kern05} thus proving that our method well computes the solid's velocity coupled with the application of gravity to the solid.
As the spatial precision increase, the curves converge to $v_y = -0.47$.
These results are quite interesting for computer graphics since, even though we are using time steps 20 times bigger (for $128\times128$) than the other method and cells 10 times bigger, we still perform a robust treatment of boundaries.

\begin{figure}[t]
\centering
\subfigure[Reference vorticity isovalues from \protect\cite{kern05}.]
{
\begin{tabular}{c}
\hspace{0.02\textwidth}\includegraphics*[angle=270,width=0.3\textwidth]{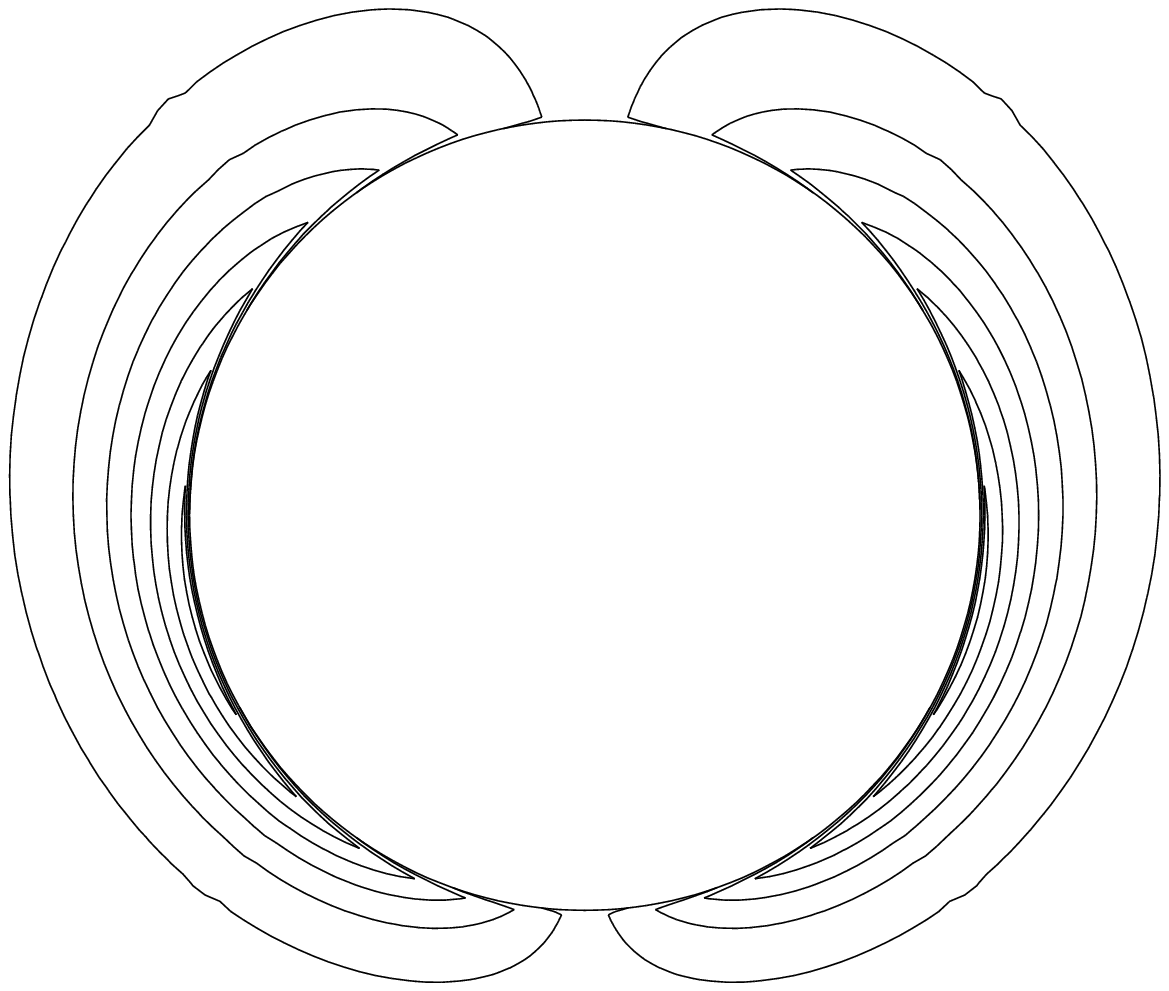}\\
\hspace{0.02\textwidth}\includegraphics*[angle=270,width=0.3\textwidth]{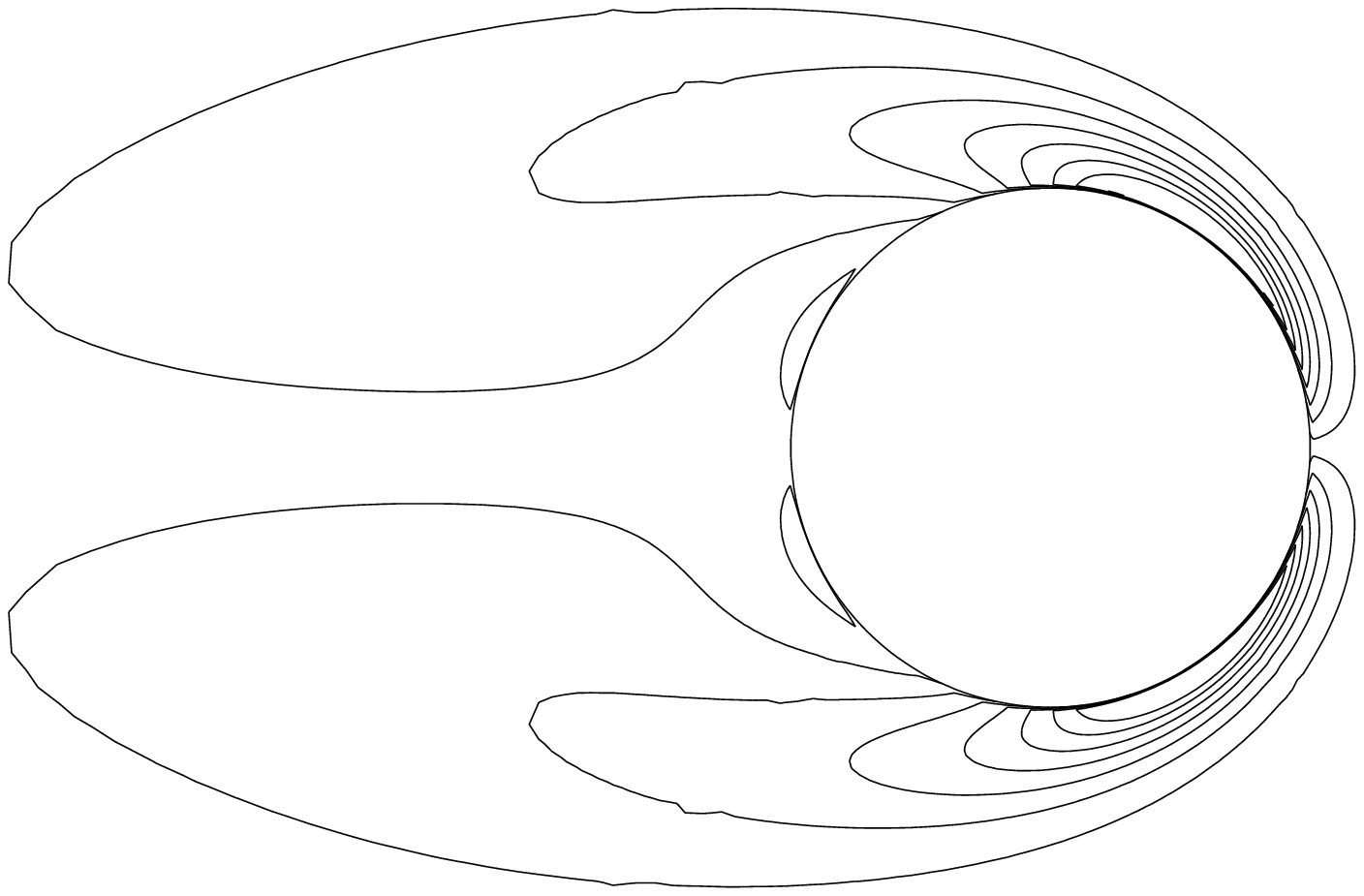}\\
\\
\end{tabular}
}\hspace{0.03\textwidth}
\subfigure[Vorticity isovalues using our method with a $300\times300$ grid.]
{\begin{tabular}{c}
\\[0.1cm]
\hspace{0.02\textwidth}\includegraphics[width=0.3\textwidth]{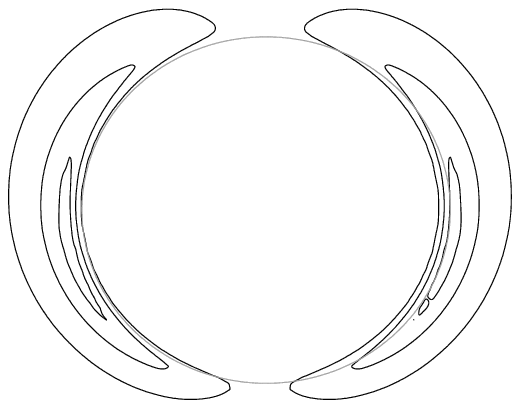}\hspace{0.02\textwidth}\\[0.15cm]
\hspace{0.02\textwidth}\includegraphics[width=0.3\textwidth]{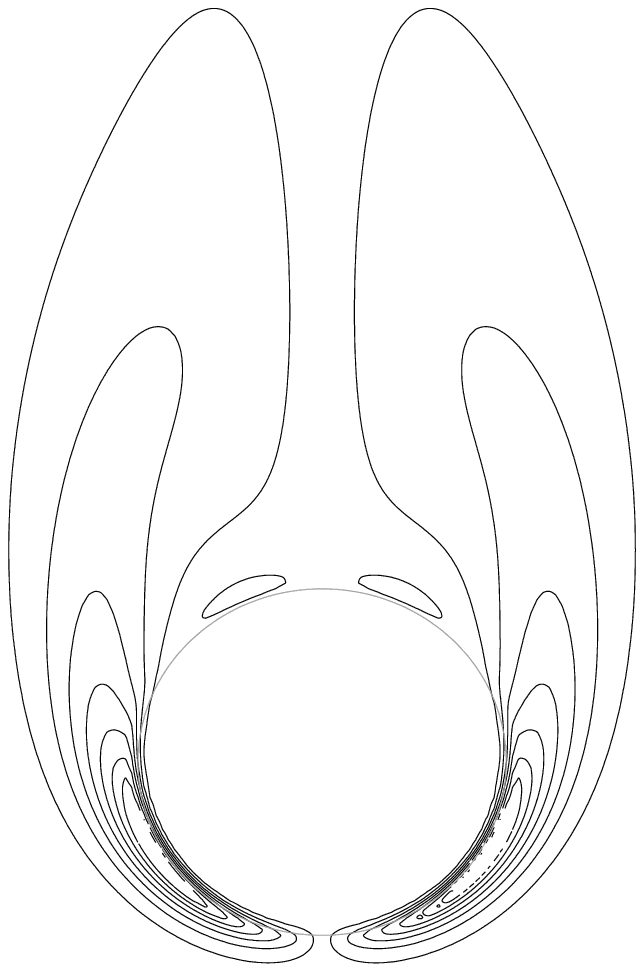}\hspace{0.02\textwidth}\\
\\
\end{tabular}
}
\caption{Vorticity isovalues at time $0.5$ (top) and $2.5$ (bottom).}
\label{fig:valid_vort}
\end{figure}

The vorticity created around the cylinder is also a revealing quantity, both physically (for fluid's dynamics) and visually (for turbulences visualisation) important.
We compare our vorticity field to the one obtained with Kern's method, isovalues of $\omega$ are shown in figure \ref{fig:valid_vort}: vorticity creation and localization are similar.
Notice that vorticity is not null inside our cylinder, this is because of the small width $\varepsilon$ used to represent the solid.
The vorticity trail created behind the cylinder is the signature of a good boundary treatment, we refer the reader to \cite{koum95} for high-resolution simulations and vorticity analysis of the flow around a cylinder with constant velocity.

\section{Implementation and results}
\label{sect:results}

We have implemented our algorithm with a Runge-Kutta 2 advection of vortex particles and second order in space differentiations.
For the following examples, we have used a low order interpolation for advecting the level set.

\begin{figure}[t]
\centering
\includegraphics[width=0.3\textwidth]{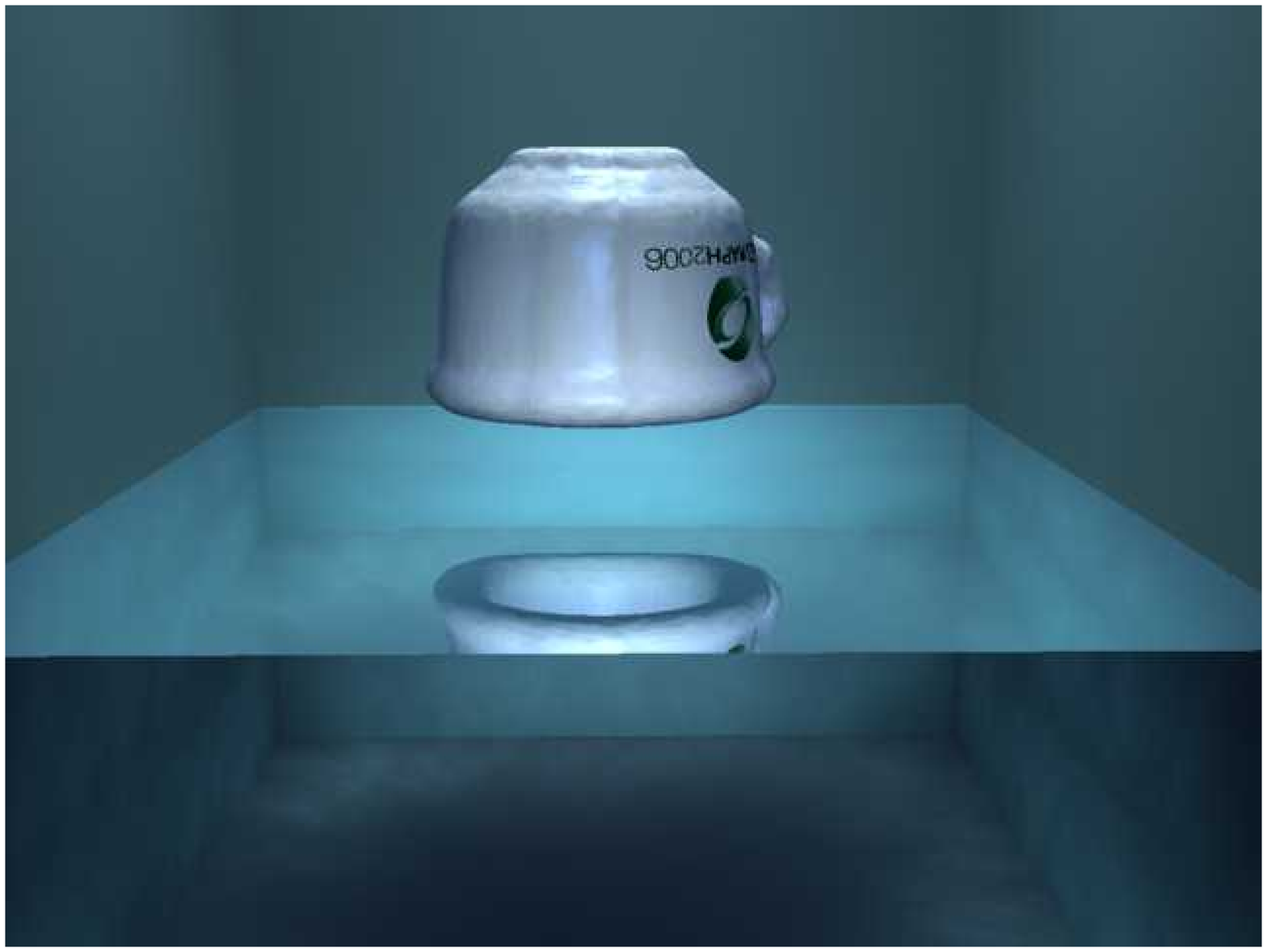}%
\includegraphics[width=0.3\textwidth]{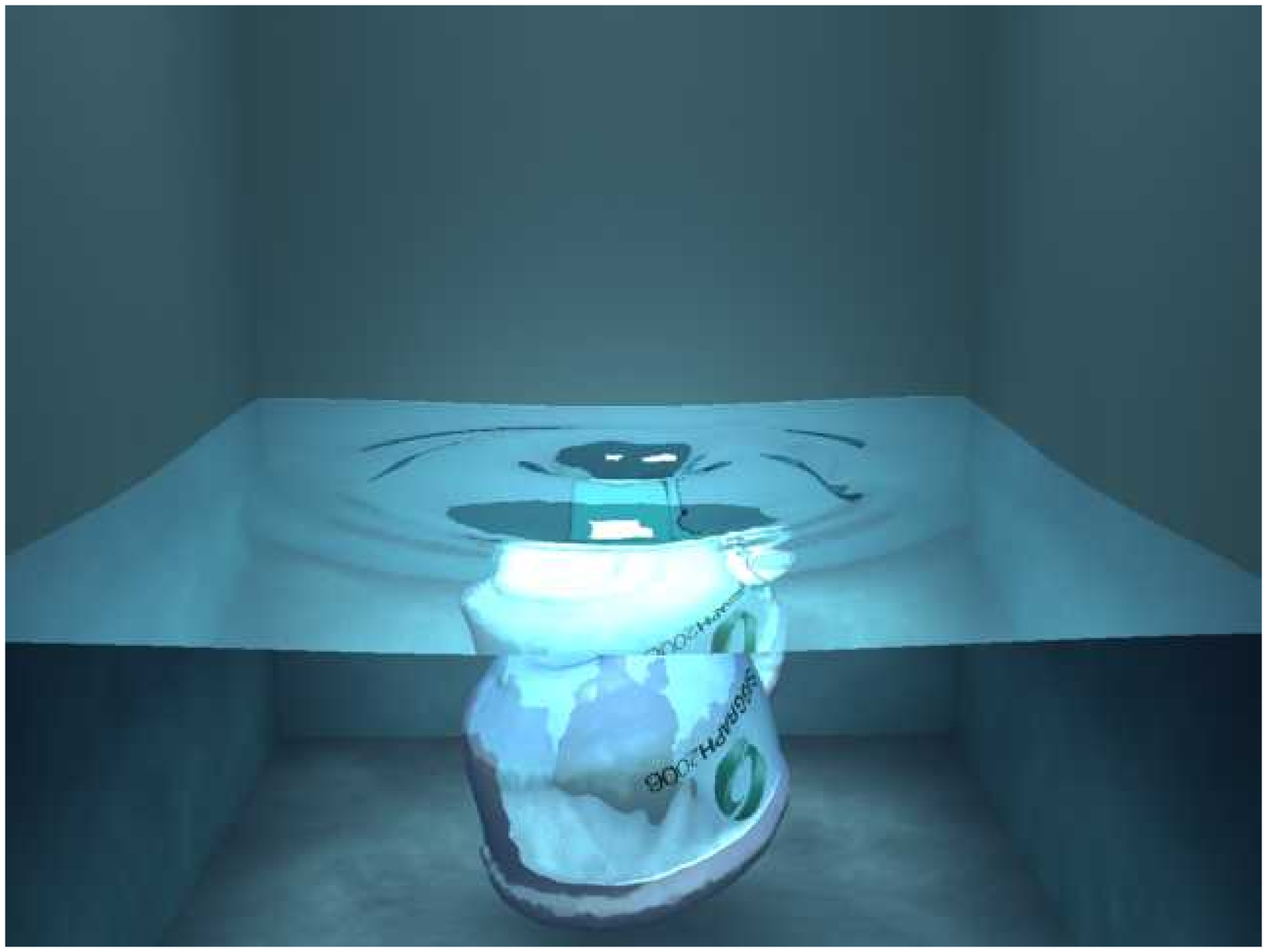}%
\\
\includegraphics[width=0.3\textwidth]{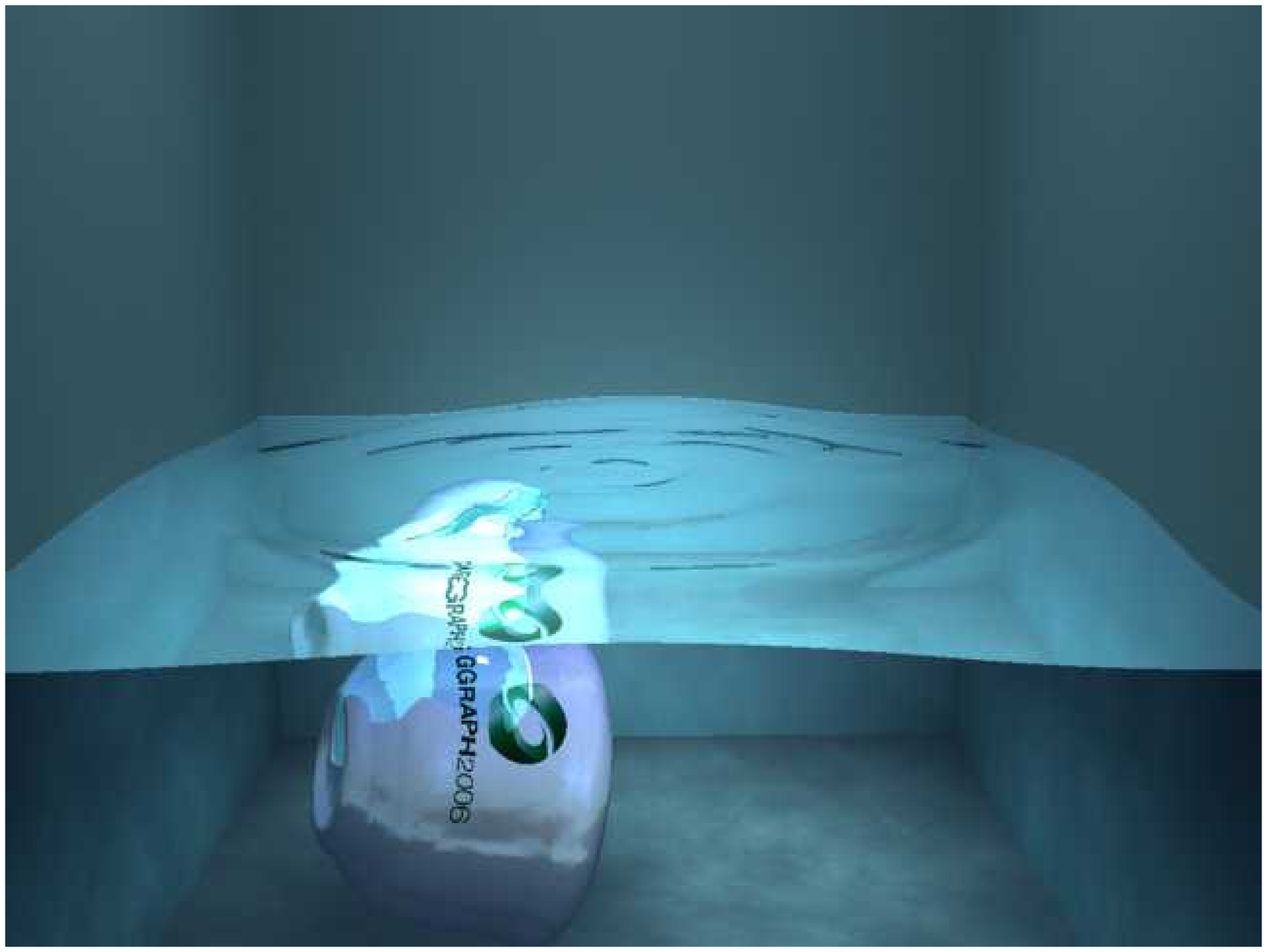}%
\includegraphics[width=0.3\textwidth]{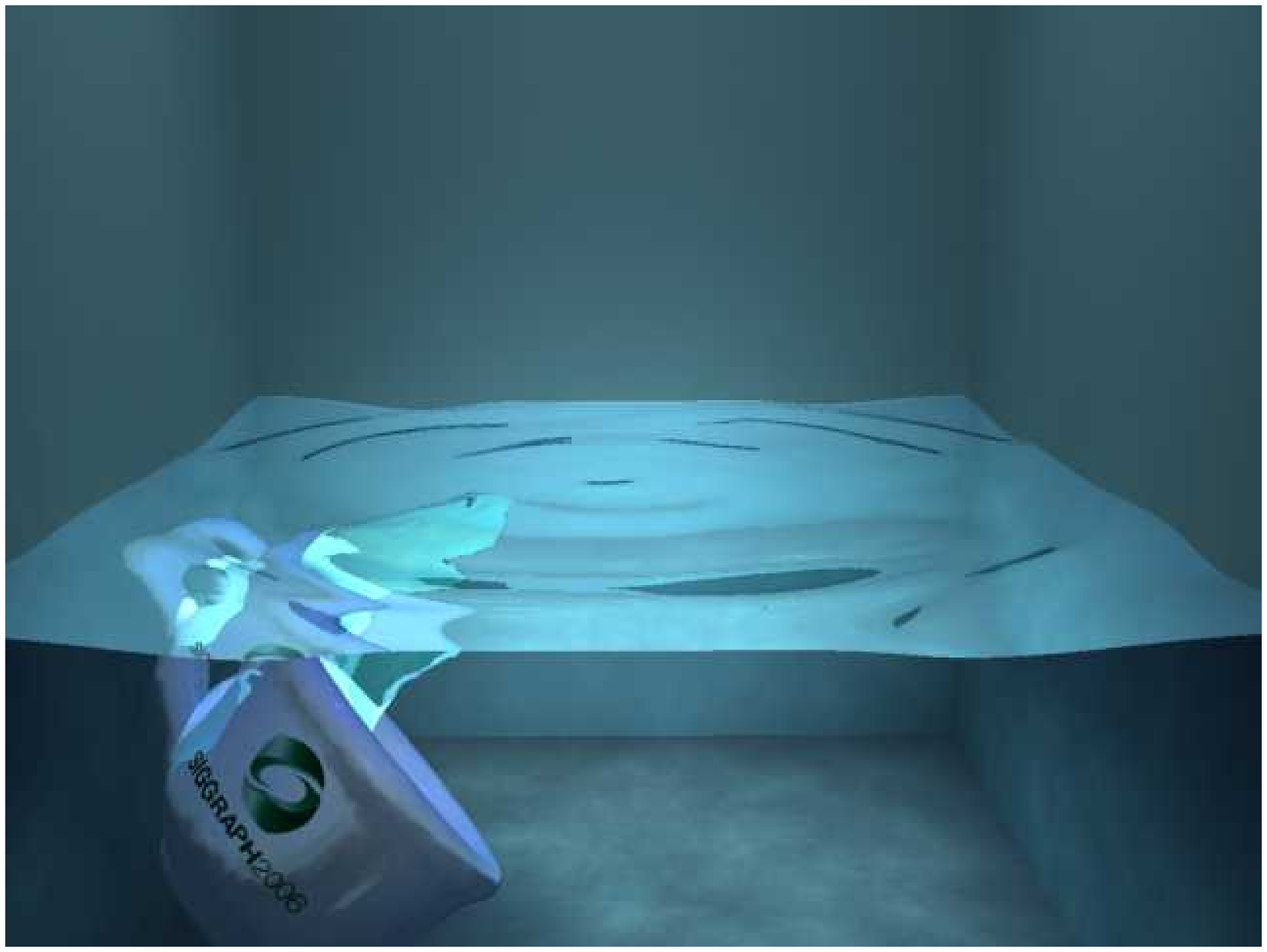}%
\includegraphics[width=0.3\textwidth]{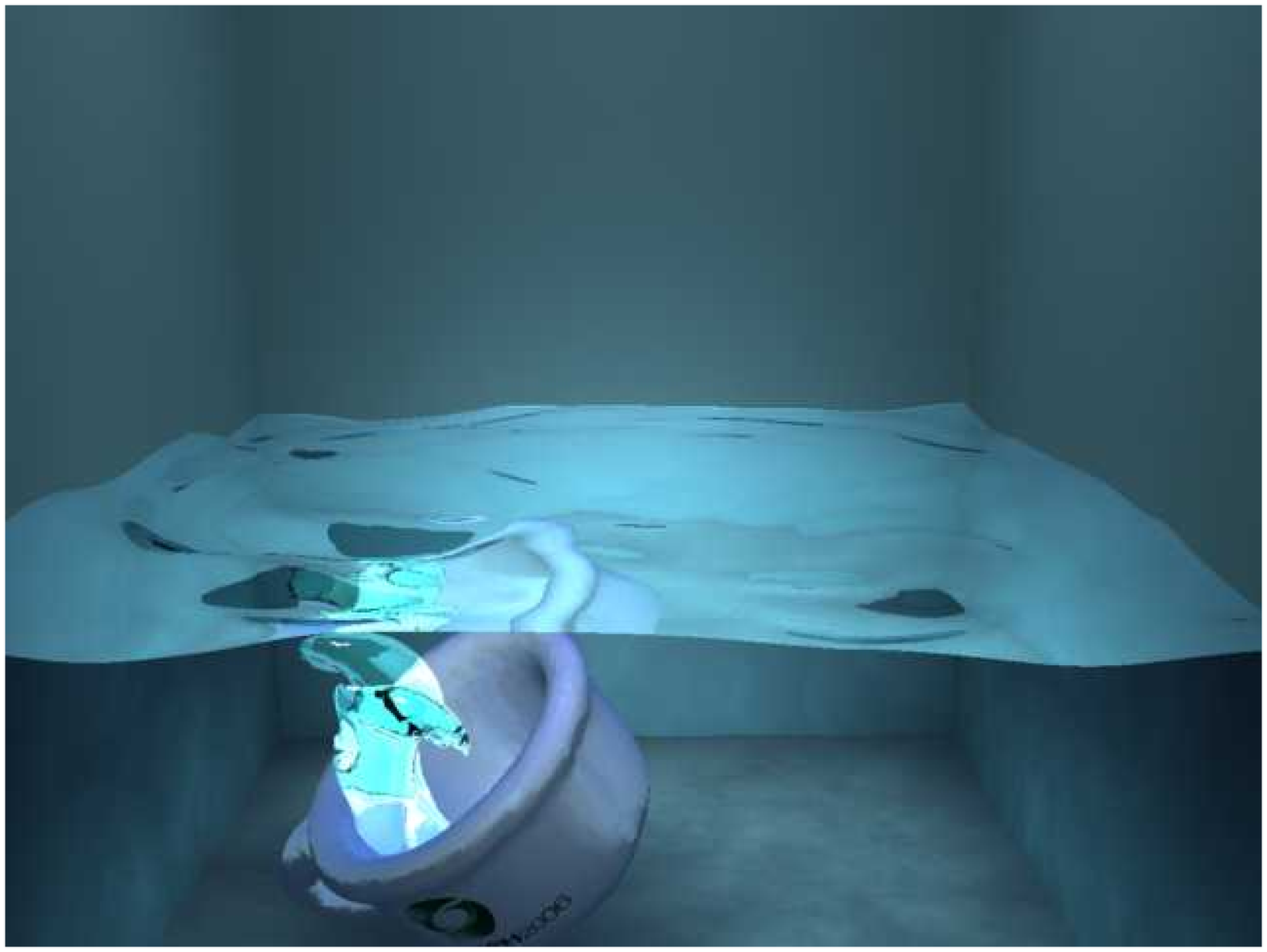}%

\caption{Interactions of a rigid object with both air and water allows
  to  simulate  complex  movements.}
\label{fig:teaser}
\end{figure}

\begin{figure}[t]
\centering
\includegraphics[width=0.63\textwidth]{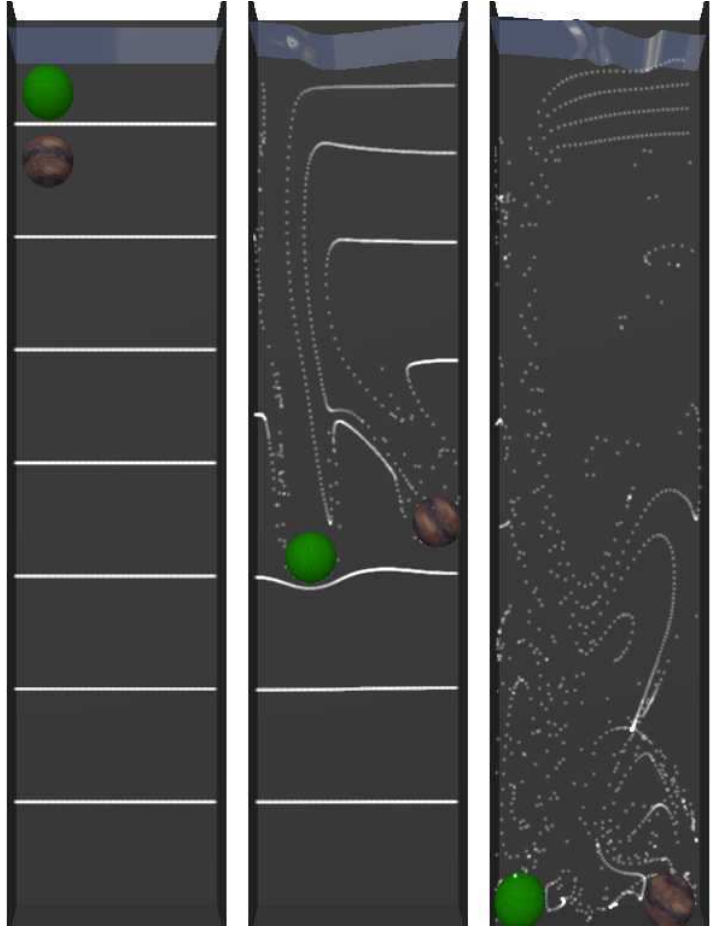}
\caption{Two sinking spheres showing drafting, kissing and tumbling effect~\protect\cite{carlson04rigid}, reproduced by our method.}
\label{fig:carlson}
\end{figure}

\begin{figure}[t]
\centering
\includegraphics[width=0.3\textwidth]{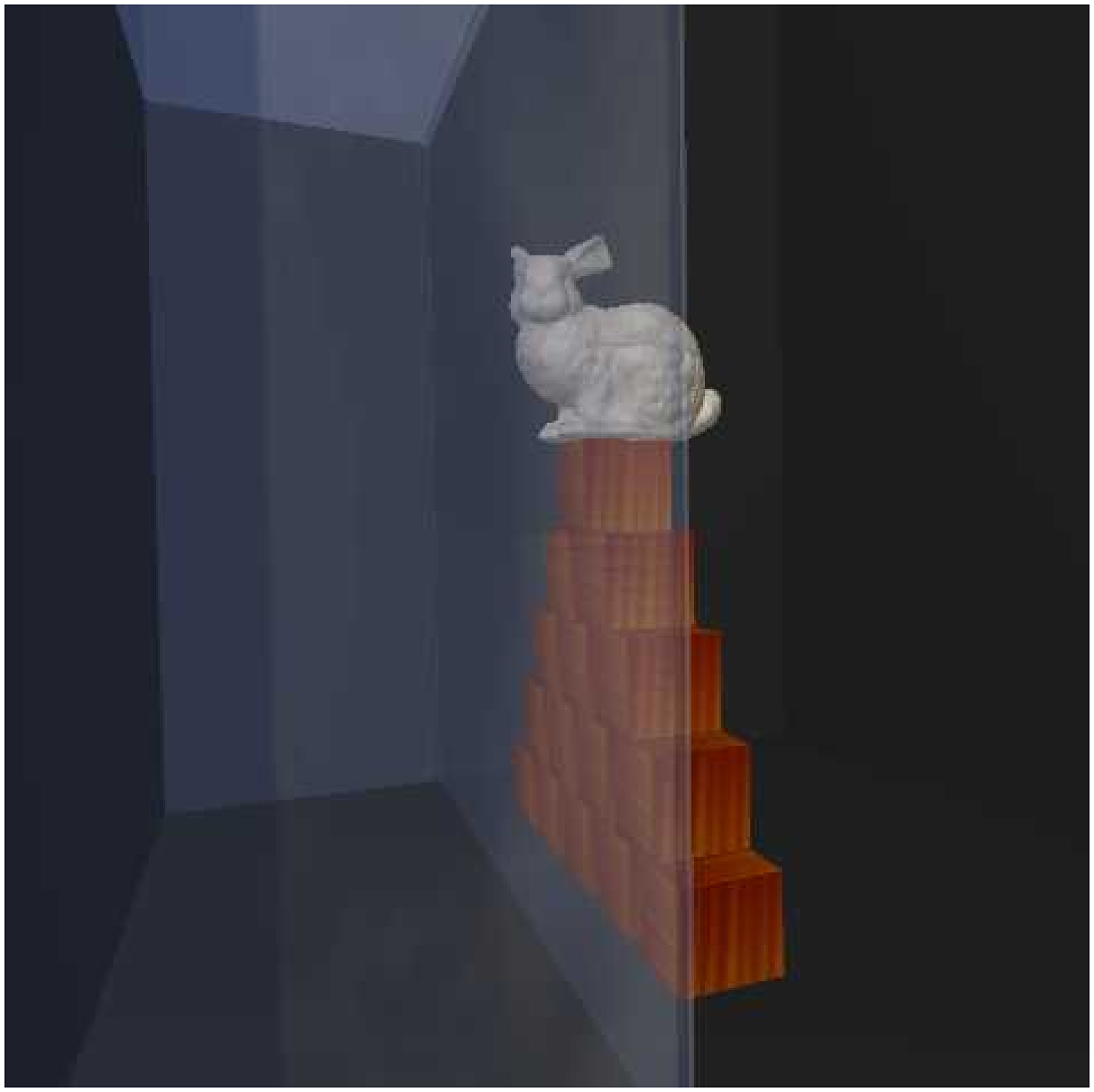}%
\includegraphics[width=0.3\textwidth]{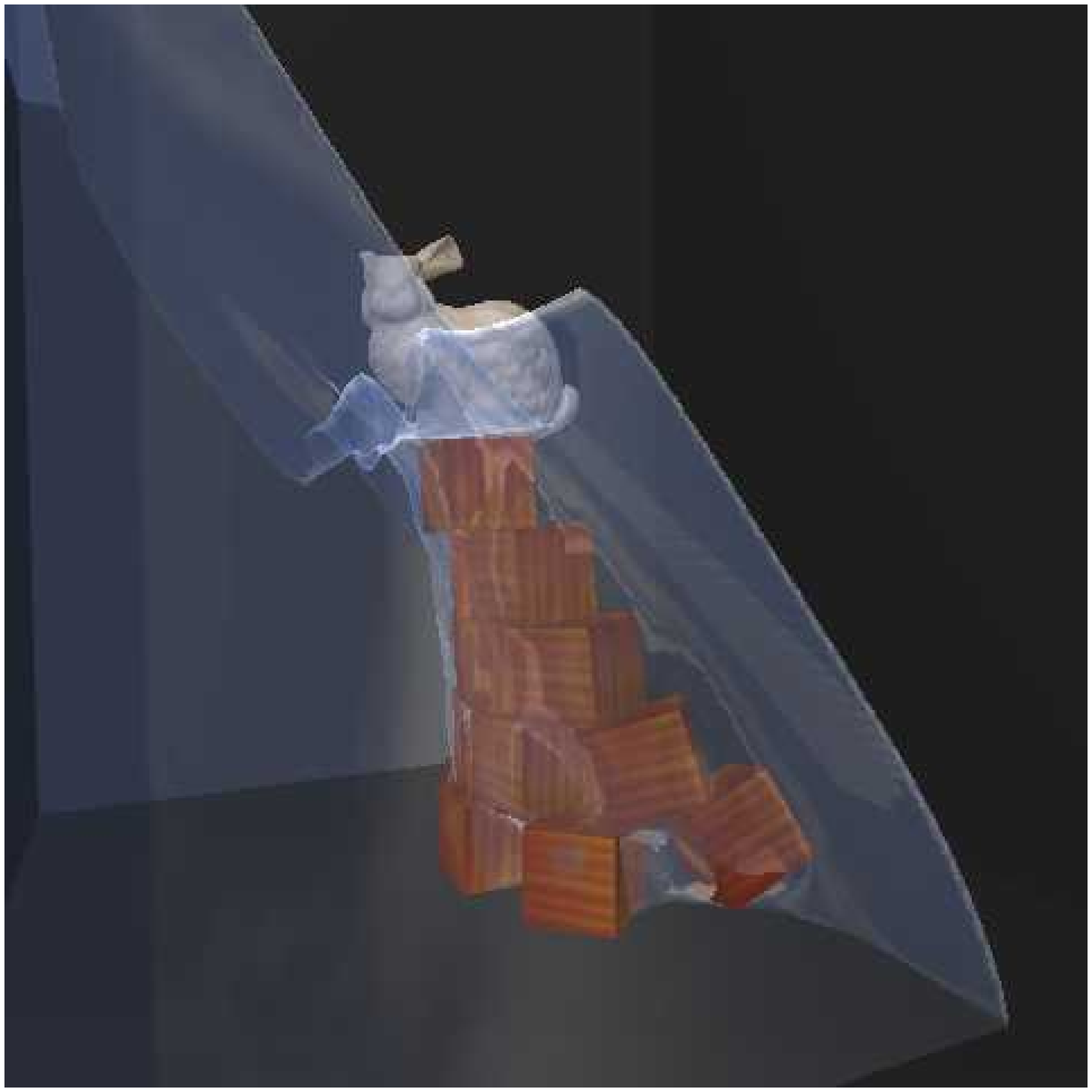}%
\includegraphics[width=0.3\textwidth]{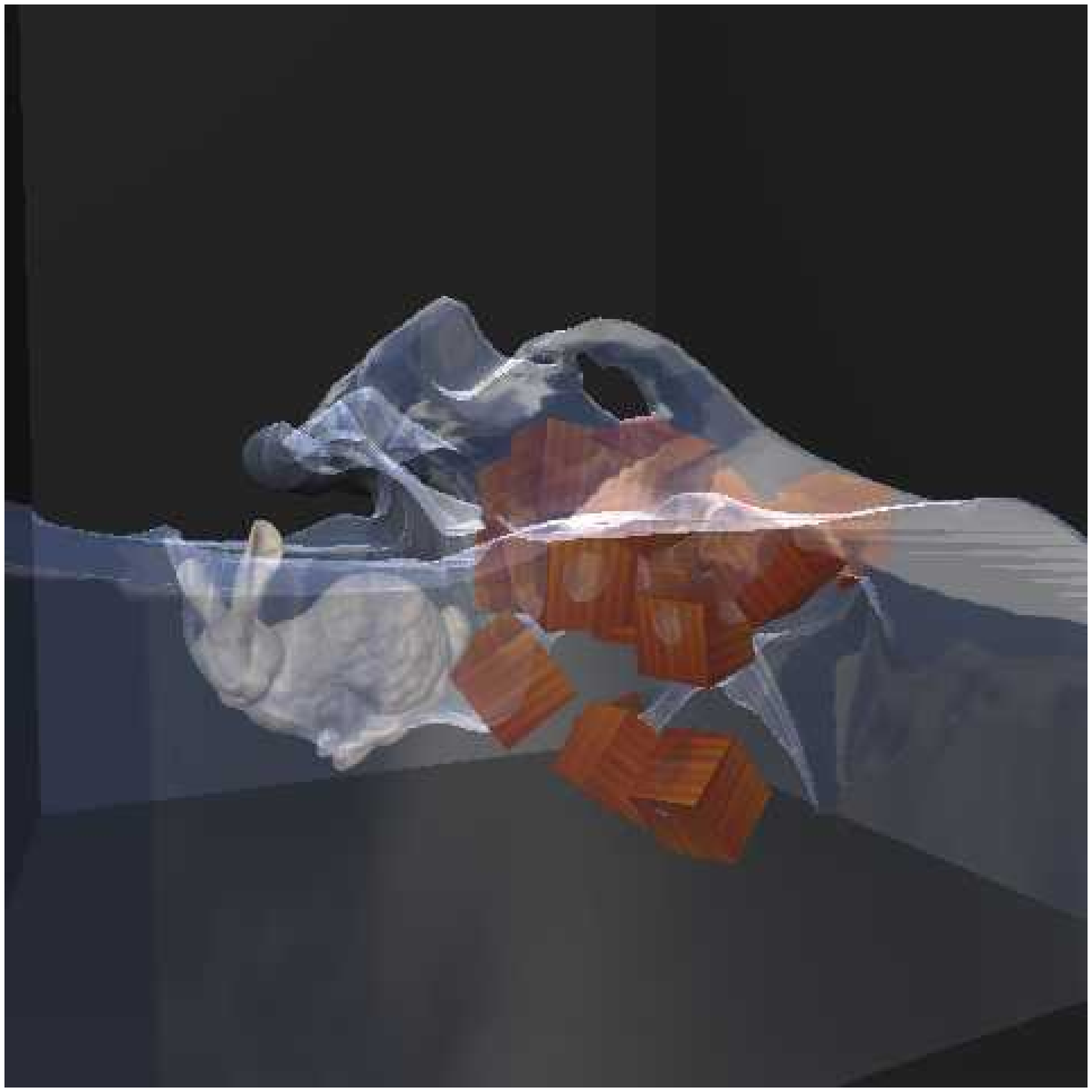}%
\caption{A pyramid encountering a wall of water.}
\label{fig:pyramid}
\end{figure}

\begin{figure}[t]
\centering
\includegraphics[width=0.23\textwidth]{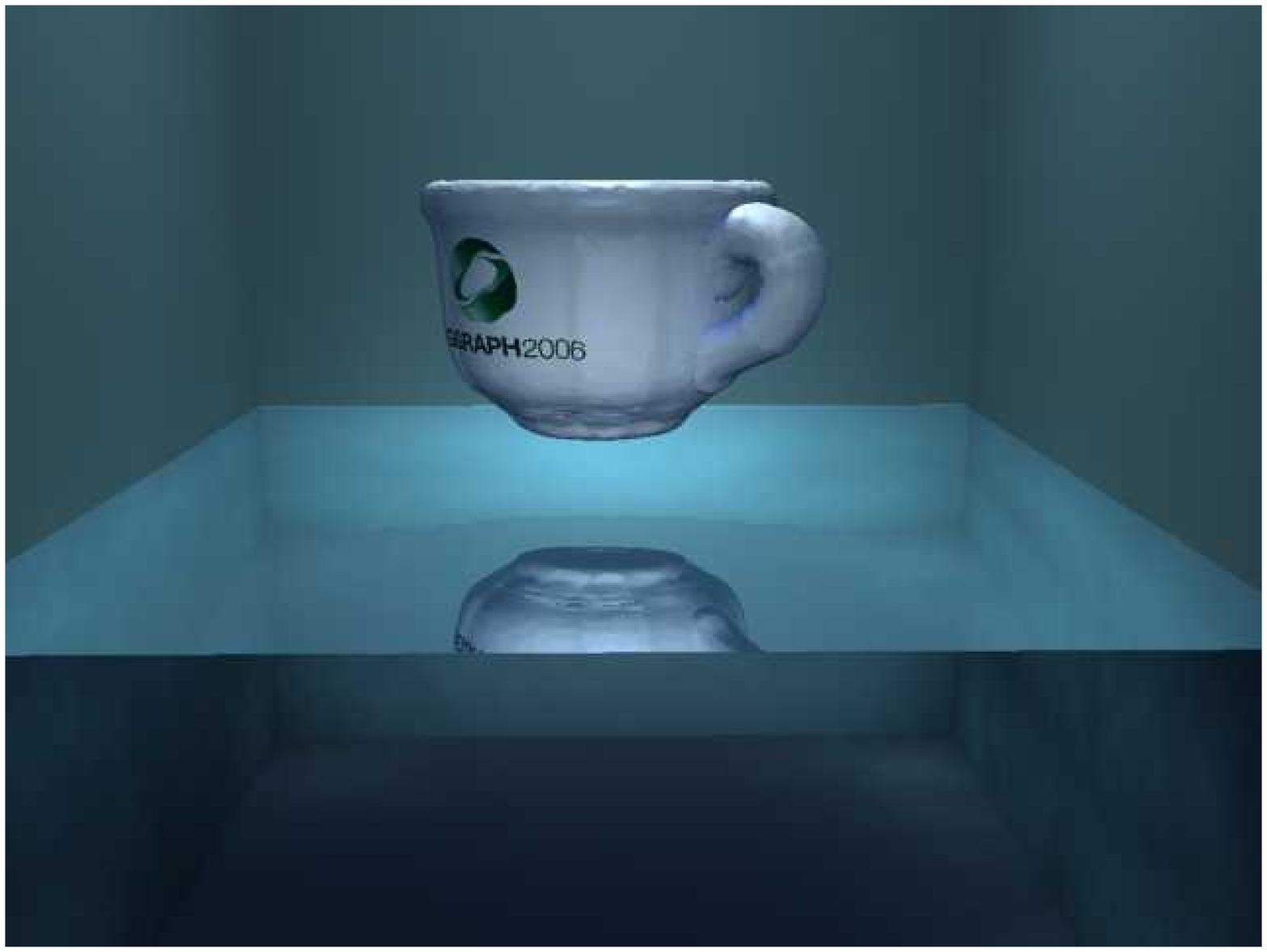}\includegraphics[width=0.23\textwidth]{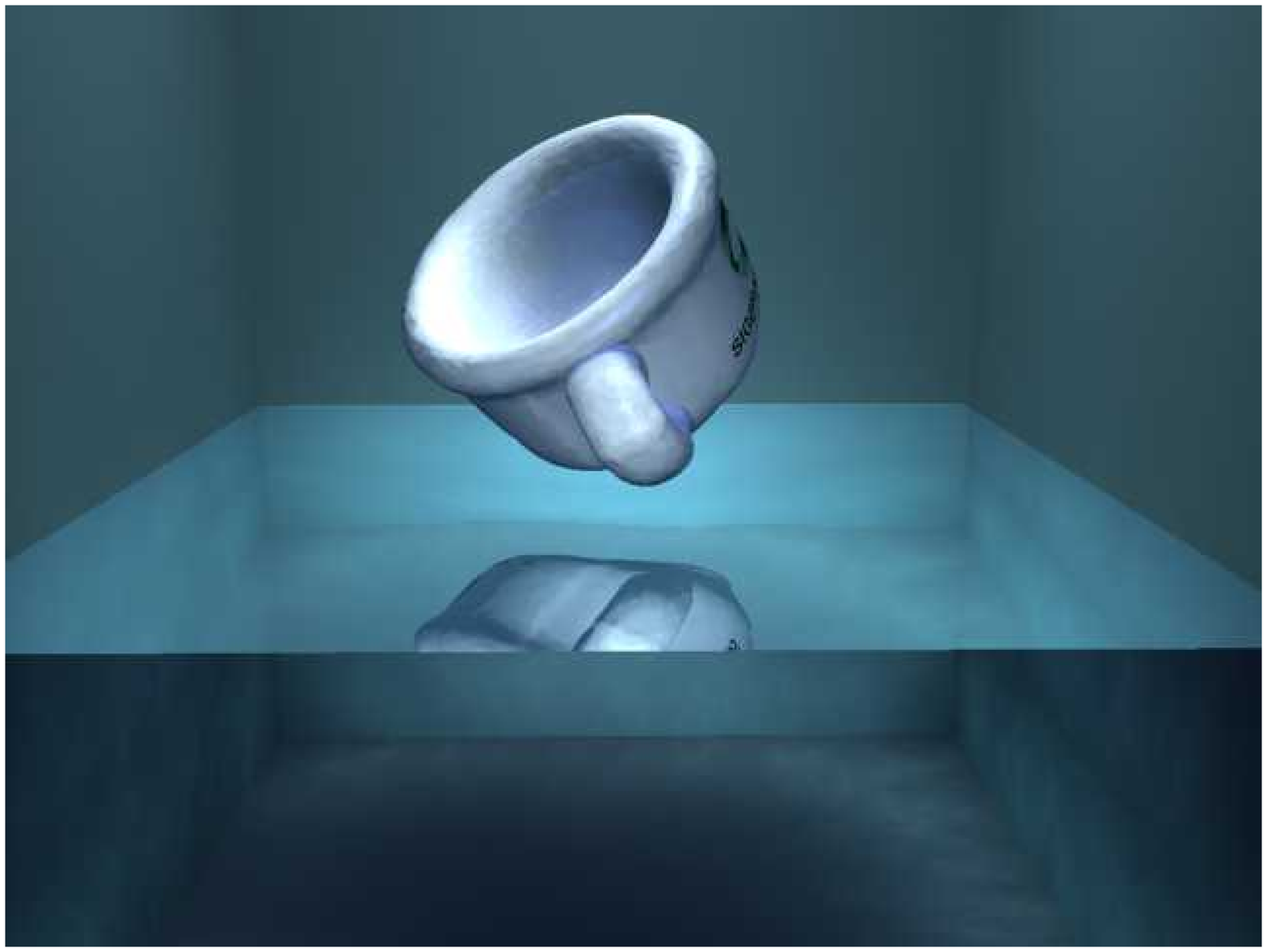}\\
\includegraphics[width=0.23\textwidth]{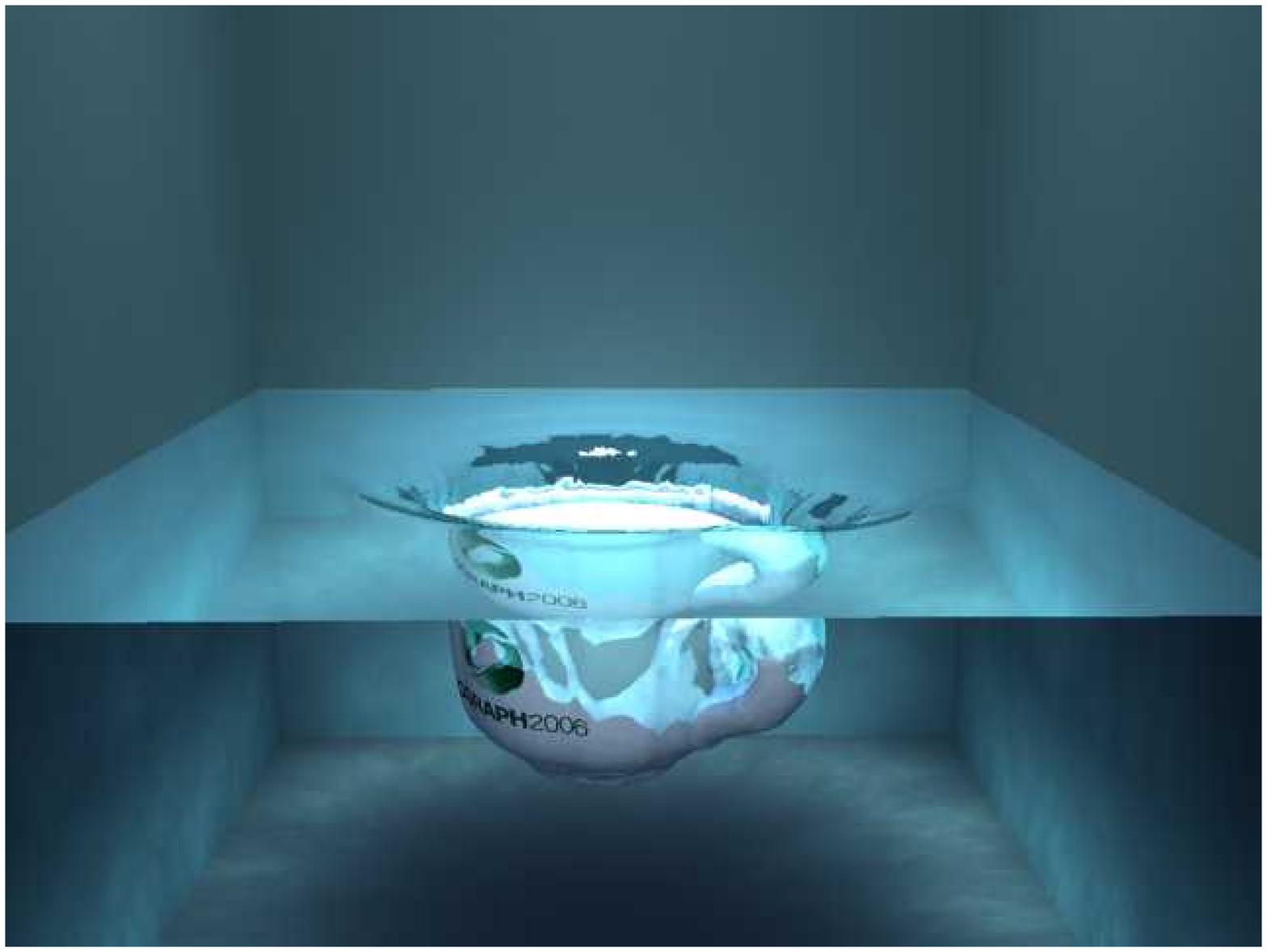}\includegraphics[width=0.23\textwidth]{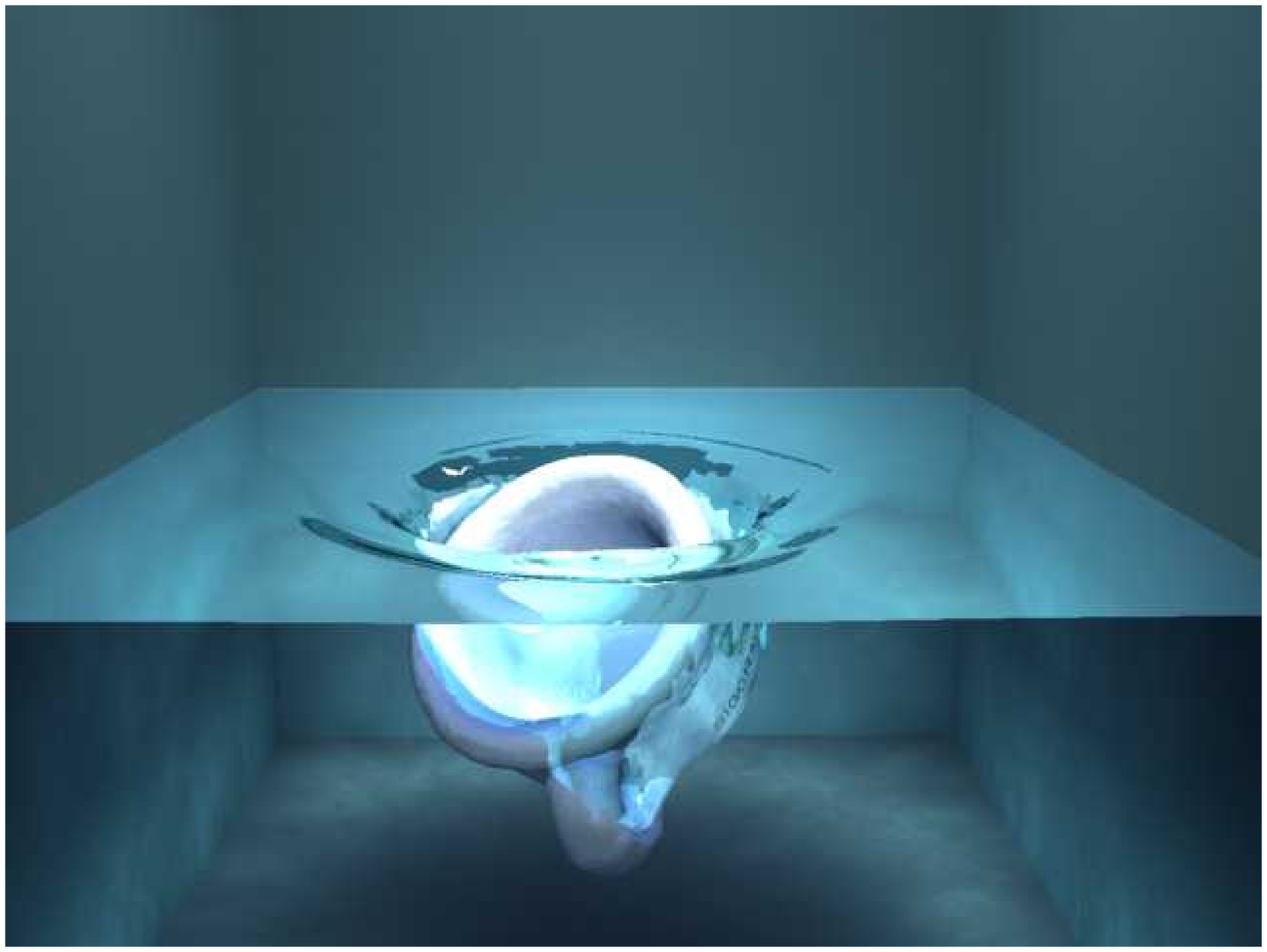}\\
\includegraphics[width=0.23\textwidth]{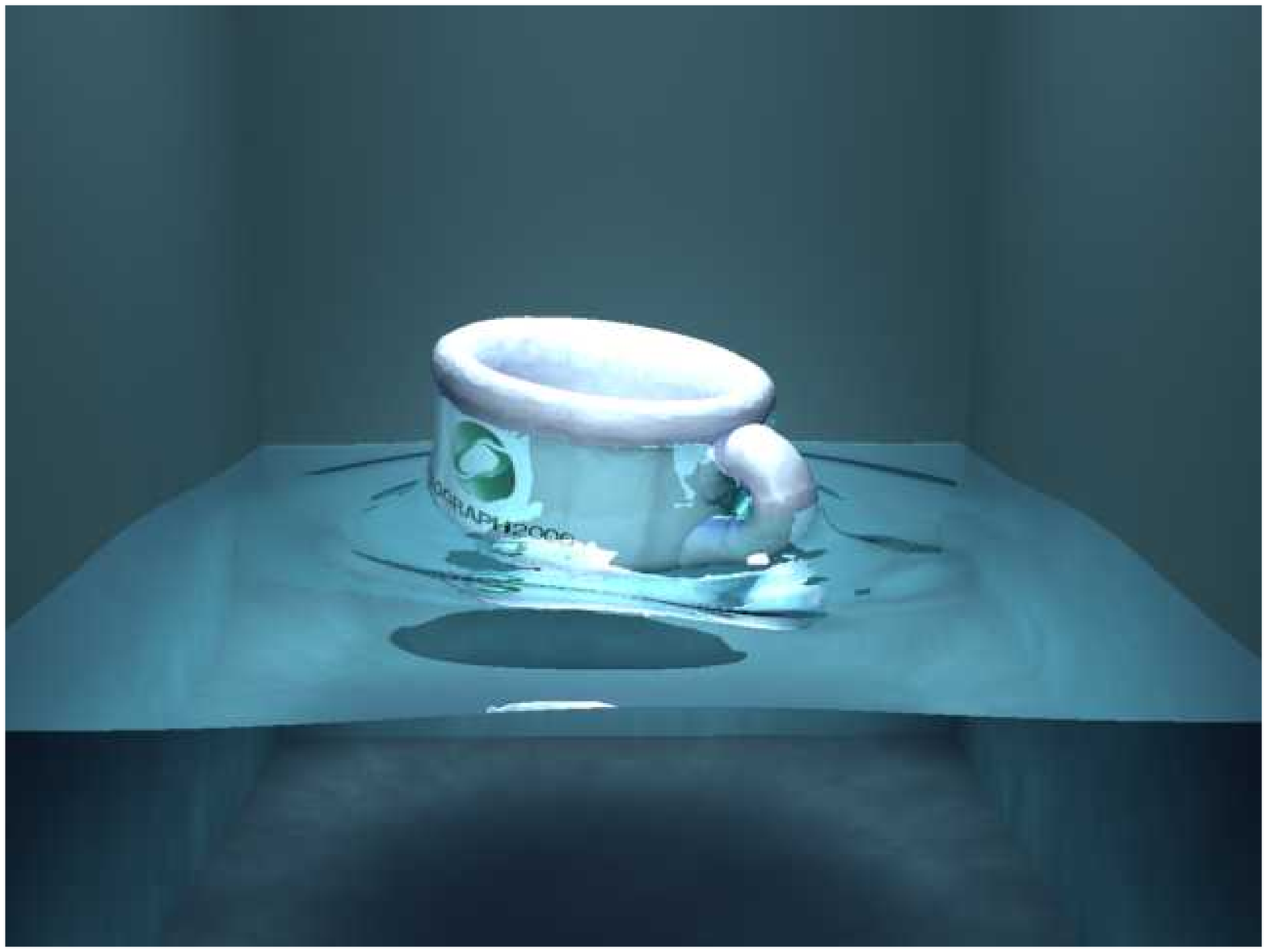}\includegraphics[width=0.23\textwidth]{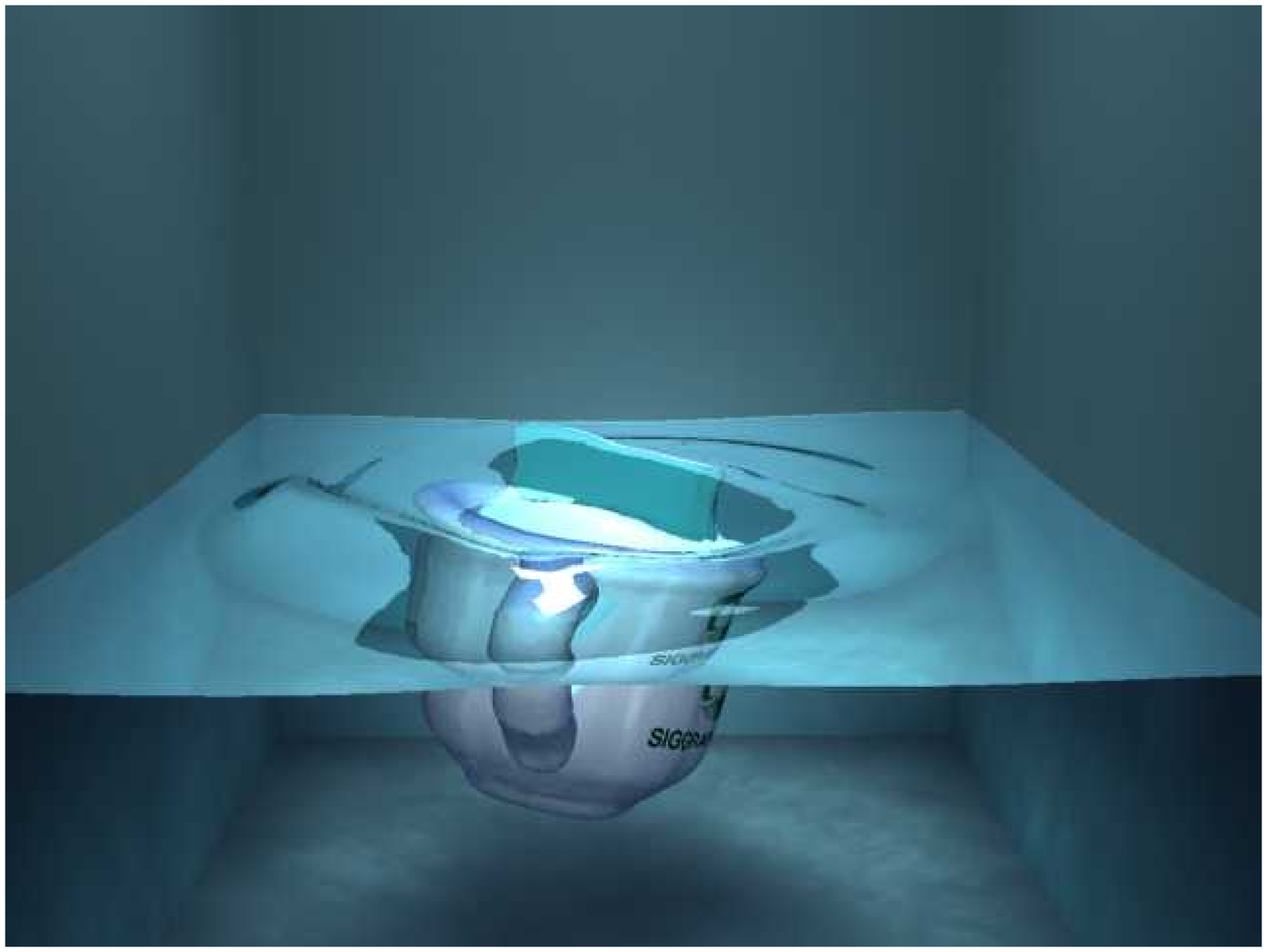}\\
\includegraphics[width=0.23\textwidth]{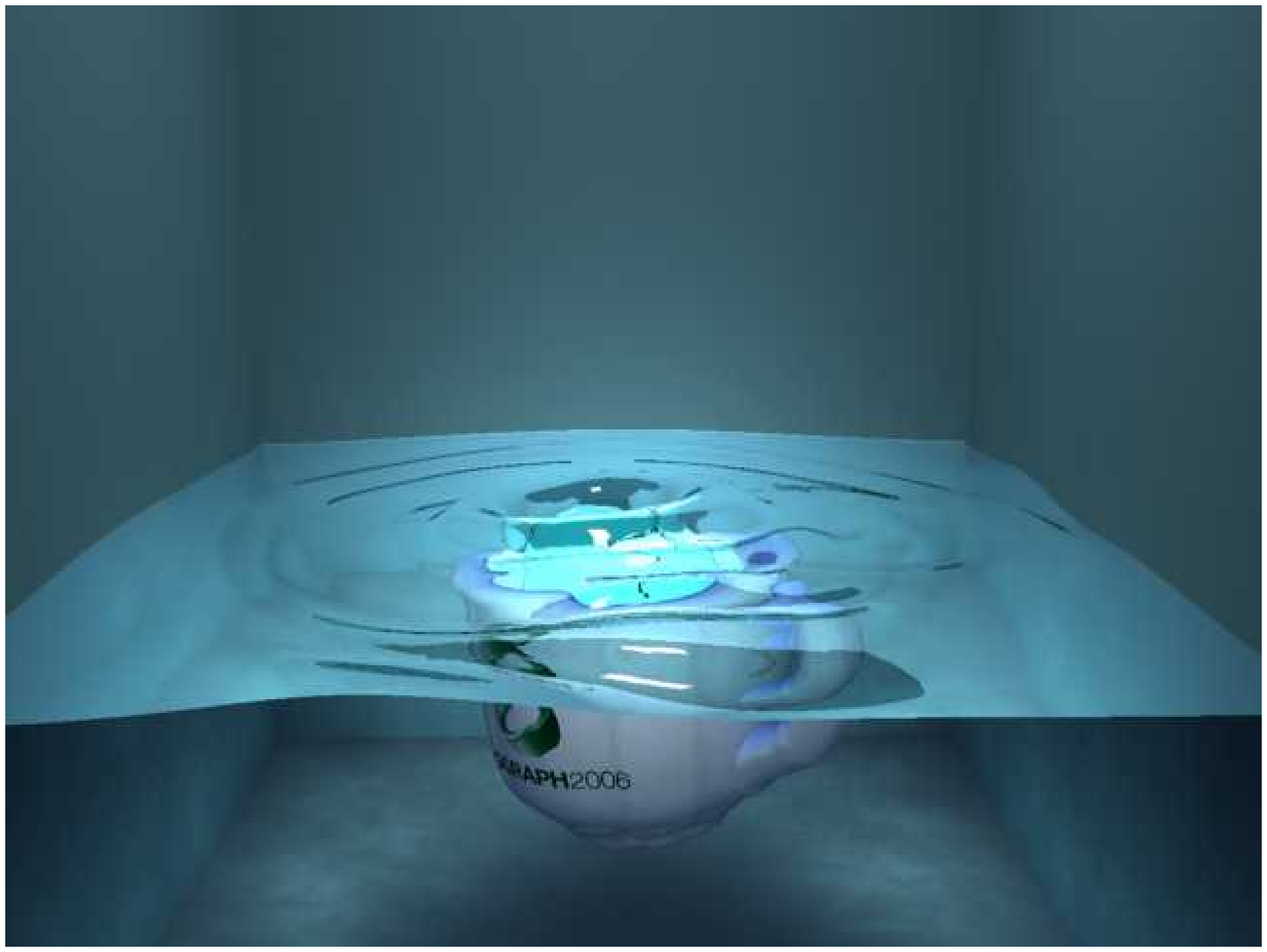}\includegraphics[width=0.23\textwidth]{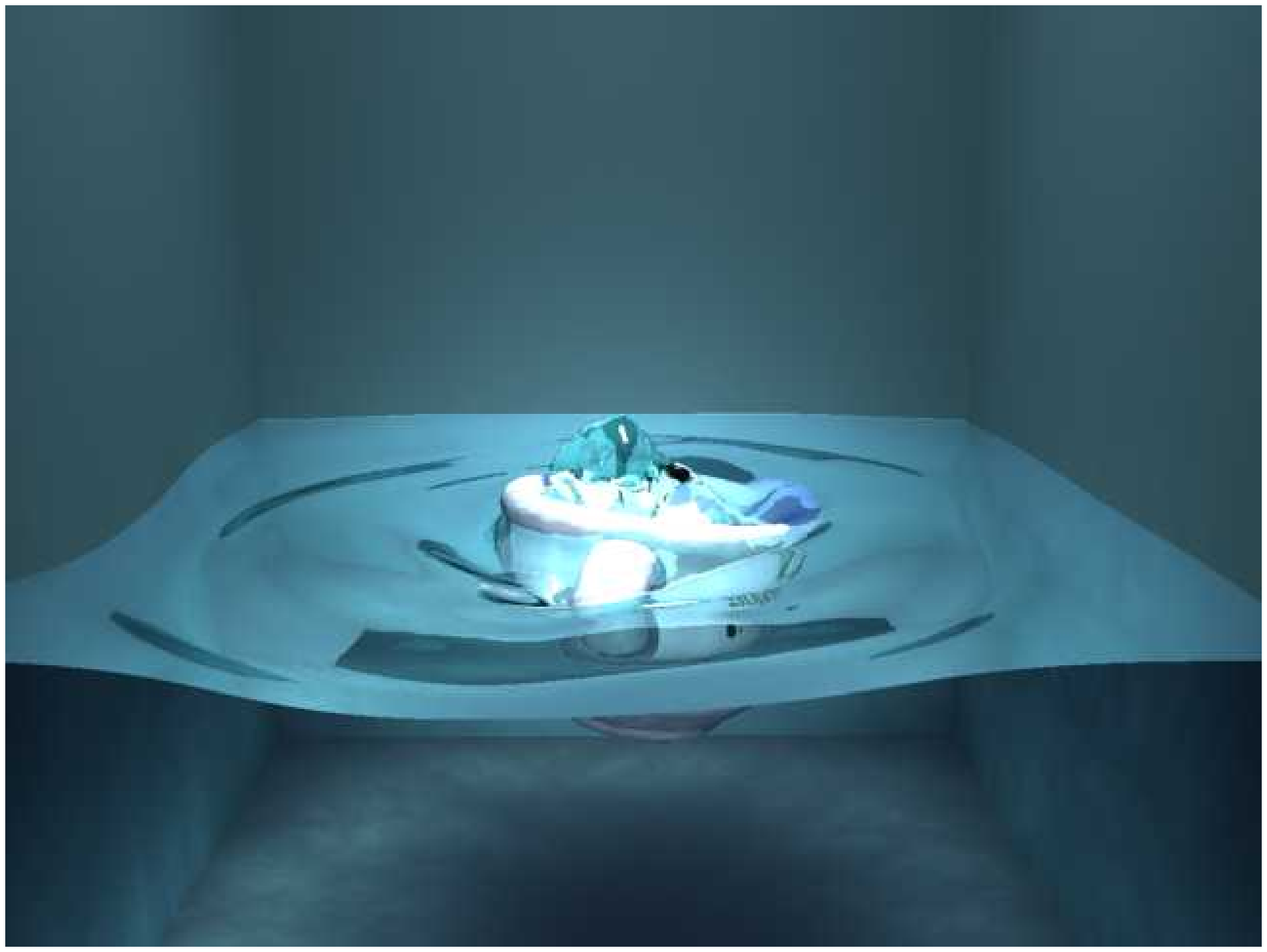}\\
\includegraphics[width=0.23\textwidth]{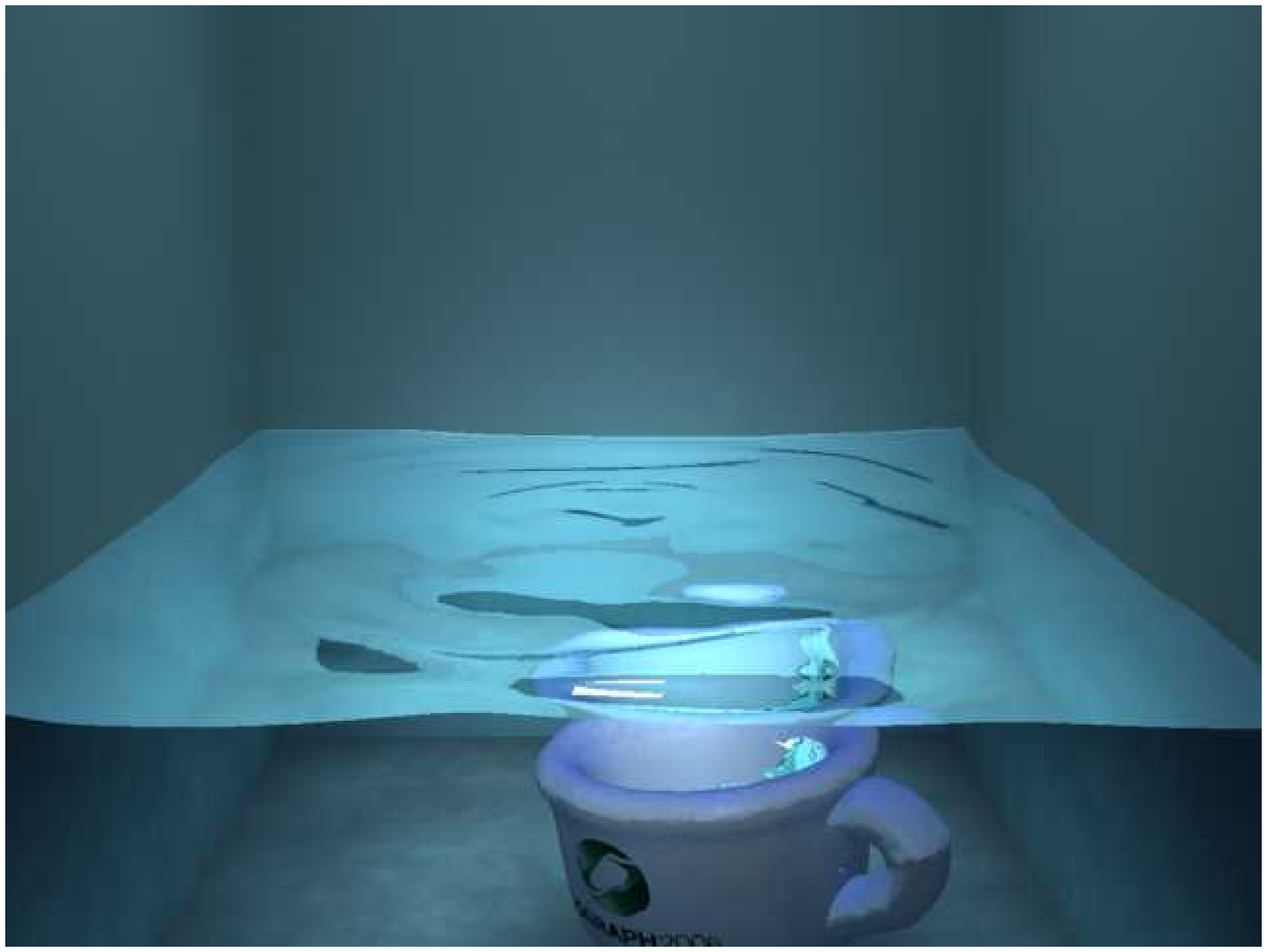}\includegraphics[width=0.23\textwidth]{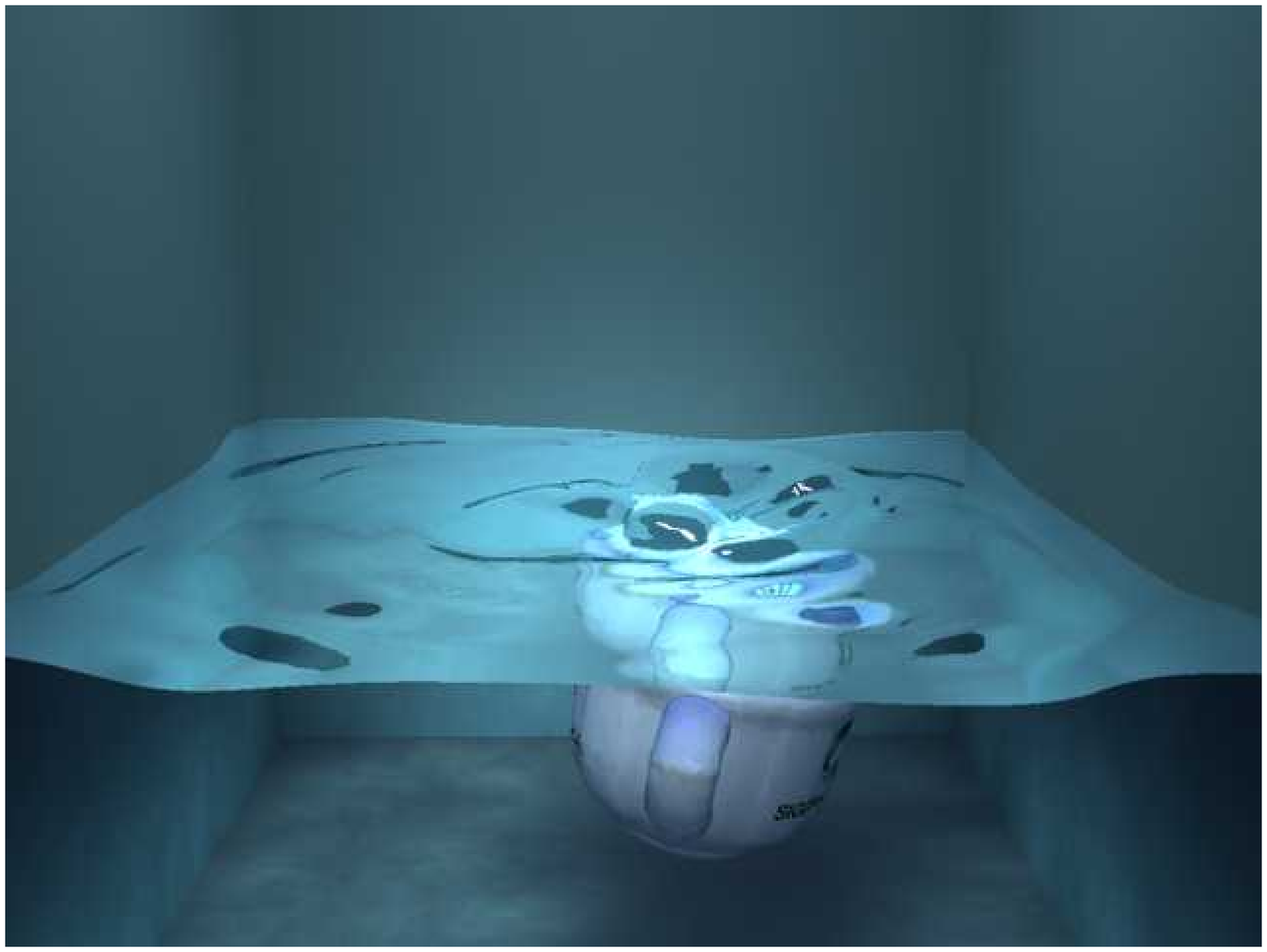}\\
\caption{Two cups enter the fluid with different orientation.} 
\label{fig:cups}
\end{figure}

In all the results we present now, we use a common gravity $\mathbf{g}=-10\mathbf{e}_z$, viscosity for water is $\nu_{water}=1.0\cdot 10^{-6}$, for air $\nu_{air}=0.82\cdot 10^{-6}$, water's density is set to $1$ and air's density to $0.001$.
The surface tension coefficient is $10^{-6}$.
All our results were computed on an Opteron 2GHz with 2Go of memory.
\\
We present different types of simulations, in figures \ref{fig:teaser} and \ref{fig:cups} cups are falling in water, rendering is done with MAYA\textsuperscript{TM}.
In figure \ref{fig:pyramid} an initialy vertical volume of water is falling on several solids and in figure \ref{fig:carlson}, we have made two spheres fall in a tank of water, similarly to~\cite{carlson04rigid}.
For those last two examples, the rendering was done with OpenGL.
A summary of the performance of the method  is provided in table \ref{table:perfs}.
\\

In the first example, shown in figure \ref{fig:teaser}, the cup has density $\rho_{cup}=1.5$.
We see that, when entering the liquid, the cup encapsulates air which makes it roll over after a while, as the captured volume tends to ride up due to gravitational forces.
One can see in the last two images of the clip, the big bubble escaping and merging back with the air.
Figure \ref{fig:cups} shows similar simulations with the same cup starting with different orientations.
One can observe the effect of surface tension and air capturing on the dynamic of the cups and fluids.
The grid's dimension is $100\times100$ and the time step used for the two simulations was $0.01$.
The total time spent for 1300 iterations per simulation was approximately 3 hours.
\\
As a second example, we took a ``wall'' of water which, under gravity, falls down and breaks the construction (see \ref{fig:pyramid}).
One can see in the third image the creation of a wave on the top of the water surface.
This wave then breaks and merges with the water under it.
\\
Figure \ref{fig:carlson} shows two spheres falling in water, this case was inspired from~\cite{carlson04rigid} and we are using the same grid resolution.
While obtaining similar dynamics for solids and fluids (the turbulent flow is observable by looking at the white ``dust'' particles), we have computed this sequence with a smaller time step of $0.01$, 400 iterations were performed in 40 minutes which represent approximately a gain of 58.
\\

One can see in the performance synthesis (table \ref{table:perfs}) that the a third of the computational cost is due to the solving of the Poisson equation, while only a small amount of time is dedicated to the particles.
The second costly part of the algorithm stands for the level set operations.
The rigid/fluid coupling takes a little more time to compute than grid-based finite differences computations such as the stretching or the viscous term computation.

\begin{table*}[tb]
\centering
\small
\begin{tabular}{|l|c|c|c|c|c|c|c|c|c|c|clcl}
\hline
Sequence & Grid Resolution         & Duration & Time       & CPU Time           \\
         &                         &          &      Steps &          / Step      \\
\hline
Cup 1     & $100\times100\times100$ & 10 s     & 1300       & 8.60 s     \\
Cup  2   & $64\times64\times64$    & 10 s     & 1000       & 2.44 s      \\
Spheres  & $68\times24\times292$   &  4 s     &  400       & 5.96 s   \\
Pyramid  & $80\times80\times80$    &  3 s     &  600       & 12.4 s   \\
\hline
\end{tabular}
\\
\vspace{0.5cm}
\begin{tabular}{|l|c|c|c|c|c|c|c|c|c|c|clcl}
\hline
Sequence &Stream       & Particules & Grid & Level-Set & Rigids          & Rigid        & Surface         & Other   \\
         &                                Solve &            &      &           &        Coupling &       Solver &         Tension &         \\
\hline
Cup   1   &  32.4 \%      & 4.93 \%    & 18.6  \% & 23.3 \% & 8.28 \%        & 3.08 \%      & 5.78 \%         & 3.66 \% \\
Cup    2  &  32.1 \%      & 4.60 \%    & 19.5  \% & 21.0 \% & 8.17 \%        & 5.38 \%      & 5.58 \%         & 3.69 \% \\
Spheres  &  42.2 \%      & 10.32 \%   & 21.18 \% & 21.2 \% & 8.35 \%        & 1.23 \%      & -               & 3.04 \% \\
Pyramid  & 33,1 \%      & 2.43 \%    & 6.81 \%  & 25.5 \% & 24.5 \%        & 1.06 \%      & 2.05 \%         & 4.50 \% \\
\hline
\end{tabular}
\caption{Performance measurements. \textit{Particles} stands for particle RK2 advection and distribution, \textit{Grid} for computation on grid (steps 3 and 7), \textit{Level-Set} for level set advection, diffusion and reinitialiation.}
\label{table:perfs}
\end{table*}

\section{Conclusion and outlook}

We have presented the first vortex particle method that handles a bi-phase fluid interacting with animated rigid bodies.
Fully based on the vorticity formulation of the Navier-Stokes equations, our method takes benefits of the advection of the vortex particles to achieve precision and robustness for large time steps.
Meanwhile, it relies on a 3D grid for efficiently computing some of the terms and redistributing the particles at each time step.
Both the bi-phase fluid and the moving solids are treated using the vorticity equations for a fluid of variable density, the 
associated level sets being advected by the fluid.
Moreover, the rigid motion of the solid is accurately computed by ensuring the continuity of velocities with the surrounding fluids.

As our results show, this method brings the two benefits of being a numerically accurate and an efficient way to simulate fluids interacting with rigid bodies.

Our level set approach for modeling the interface between the fluids and immersed bodies could be extended to other cases. One can for instance enrich the type of physics underlying the fluid-body interaction, while keeping the simplicity provided by the immersed interface approach and its particle discretization. Using the framework developed in \cite{CotMai06} it is actually possible to handle flexible bodies interacting with fluids, a possibility that we plan to explore.
Another idea for future work is to implement multi-level particle methods in our fluid-body models. One may envision two types of multi-level discretization. One would be to capture the interfaces, that is the level-set functions, and the flow, that is the vorticity  with different grid resolution. Using finer resolution for the level-set function should allow to capture finer details on the interfaces with little computational overhead. Indeed it should be emphasized that, since the time-step of particle methods is only constrained by the amount of strain in the flow, refining the resolution for the level-set functions does not result in a smaller time-step. Alternatively, one may consider using a full (for the flow and the interfaces) multi-level vortex method, in the sprit of AMR methods, along the lines of \cite{BerCotKou05}. Fluid-body interactions is a challenging problem  to clarify the respective benefits of these two approaches.

\bibliographystyle{acm}
\bibliography{article}

   \end{document}